\documentclass{scrartcl}

\usepackage{xcolor}

\usepackage{amsmath}
\usepackage{amsthm}
\usepackage{amssymb}
\usepackage{amsfonts}
\usepackage{mathrsfs}
\usepackage{color}
\usepackage[utf8]{inputenc}
\usepackage{enumitem}
\usepackage{multicol}
\usepackage[english]{babel}
\usepackage[babel]{csquotes}
\usepackage{mathtools}
\usepackage{rotating}
\usepackage{latexsym}
\usepackage{epsfig}
\usepackage{verbatim}
\usepackage{thmtools,thm-restate}
\usepackage{listings} 

\usepackage{makecell,booktabs,amsxtra}

\usepackage{adjustbox}
\usepackage{ifthen}
\usepackage[position=top]{subfig}
\usepackage{emptypage}
\usepackage{pgfplots}
\pgfplotsset{xticklabels={}, yticklabels={} width=7cm,compat=1.8}
\usepackage{pgfplotstable}
\pgfmathsetseed{1138} 
\pgfplotstableset{ 
	create on use/x/.style={create col/expr={42+2*\pgfplotstablerow}},
	create on use/y/.style={create col/expr={(0.6*\thisrow{x}+130)+5*rand}}
}
\pgfplotstablenew[columns={x,y}]{30}\loadedtable

\usepackage[ruled]{algorithm2e}

\SetAlFnt{\small}
\SetAlCapFnt{\small}
\SetAlCapNameFnt{\small}
\SetAlCapHSkip{0pt}
\IncMargin{-\parindent}

\usepackage{url}
\usepackage[square,numbers,sort]{natbib}

\usepackage{pgfplots,pgfplotstable}
\pgfplotsset{compat=1.8}

\usepackage[bookmarks=false,colorlinks=true,linkcolor=blue,urlcolor=blue,citecolor=red,pdfborder={0 0 0.5}]{hyperref}

\theoremstyle{plain}
\newtheorem{theorem}{Theorem}[section]
\newtheorem{lemma}[theorem]{Lemma}
\newtheorem{proposition}[theorem]{Proposition}
\newtheorem{corollary}[theorem]{Corollary}

\declaretheoremstyle[
notefont=\bfseries, notebraces={}{},
bodyfont=\normalfont\itshape,
headformat=\NAME \NOTE
]{nopar}
\declaretheorem[style=nopar,name=Corollary]{corollary*}
\declaretheorem[style=nopar,name=Theorem]{theorem*}

\theoremstyle{definition}
\newtheorem{definition}[theorem]{Definition}
\newtheorem{remark}[theorem]{Remark}

\newcount \mycount
\newcount \mycountb

\DeclareMathOperator{\csp}{CSP}

\DeclareMathOperator{\A}{{\mathbb A}}

\DeclareMathOperator{\F}{{\mathbb F}}

\DeclareMathOperator{\Idemp}{{\operatorname{\mathbb Idem}}}

\DeclareMathOperator{\Pol}{Pol}
\DeclareMathOperator{\Ind}{Ind}

\DeclareMathOperator{\Maltsev}{{\operatorname{Malt}}}

\newcommand{\TS}[1]{{\operatorname{TS}(#1)}}

\providecommand{\dotdiv}{
  \mathbin{
    \vphantom{+}
    \text{
      \mathsurround=0pt 
      \protect\ooalign{
        \noalign{\kern-.45ex}
        \hidewidth$\smash{\cdot}$\hidewidth\cr 
        \noalign{\kern.45ex}
        $-$\cr 
      }%
    }%
  }%
}
\providecommand{\cupdot}{
  \mathbin{
    \vphantom{+}
    \text{
      \mathsurround=0pt 
      \protect\ooalign{
        \noalign{\kern-.4ex}
        \hidewidth$\smash{\cdot}$\hidewidth\cr 
        \noalign{\kern.4ex}
        $\cup$\cr 
      }%
    }%
  }%
}

\newcommand{\ppleq}{\mathrel{\leq}}

\theoremstyle{definition}
\newtheorem{examplex}[theorem]{Example}
\newenvironment{example}
{\pushQED{\qed}\examplex}
{\popQED\endexamplex}

\theoremstyle{definition}

\makeatletter
\newcommand\footnoteref[1]{\protected@xdef\@thefnmark{\ref{#1}}\@footnotemark}
\makeatother

\makeatletter
\newcommand*{\rom}[1]{\expandafter\@slowromancap\romannumeral #1@}
\makeatother

\makeatletter
\def\blfootnote{\xdef\@thefnmark{}\@footnotetext}
\makeatother

\usepackage{setspace}

\usepackage{epigraph}
\newcommand{\prog}[1]{\mathcal #1}

\usepackage{siunitx}
\newcounter{anzahl}
\def\nodeDist{0.4}

\newcommand{\tikzpathFreeVariables}[1]{
\begin{tikzpicture}[scale=0.5]
\setcounter{anzahl}{0}

\foreach \y/\c [count=\xi from 1] in #1 {
    \ifthenelse{\c=0}{
    \node[var-f] (\xi) at (\nodeDist*\xi,\nodeDist*\y) {};}{
    \node[var-b] (\xi) at (\nodeDist*\xi,\nodeDist*\y) {};}
    
    \setcounter{anzahl}{\xi}
    }
    
\foreach \y [count=\yi from 1,count=\yii from 2] in {2,...,\arabic{anzahl}} 
    \path (\yi) edge (\yii);
\end{tikzpicture}
}

\newcommand{\toEdge}{\mathbin
{\begin{tikzpicture}[baseline =-1mm]
    \path[->] (0,0) edge (0.37,0);
\end{tikzpicture}}}
\newcommand{\fromEdge}{\mathbin
{\begin{tikzpicture}[baseline =-1mm]
    \path[<-] (0,0) edge (0.37,0);
\end{tikzpicture}}}
\tikzstyle{var-b}=[circle,fill,draw=white,inner sep=0pt,minimum size=3.5pt]
\tikzstyle{bullet}=[circle,fill,draw=white,inner sep=0pt,minimum size=3.2pt]
\tikzstyle{var-f}=[circle,draw,inner sep=0pt,minimum size=3pt]

\usetikzlibrary{arrows,arrows.meta}
\usetikzlibrary{bending,calc}

\newcommand{\albertTwoElementPoset}{
\begin{tikzpicture}[scale=0.85]
\node[var-b,label=right:$\bC_1$] (C1) at (0,0) {};
\node[var-b,label=right:$\bP_2$] (P2) at (0,-1) {};

\node[var-b,below right  of=P2,label=right:$\bB_2$] (B2) {};
\node[var-b,below of=B2,label=right:$\bB_3$] (B3) {};
\node[below of=B3,label=right:] (Bn) {$\vdots$};
\node[var-b,below of=Bn,label=right:$\bB_{\infty}$] (Binf) {};

\node[var-b,below left of=B2,label=left:$\bD_2$] (D2) {};
\node[var-b,below of=D2,label=left:{}] (D3) {};
\node[below of=D3,label=right:{}] (Dn) {$\vdots$};
\node[var-b,below of=Dn,label=left:$\bB_{\infty}^{\leq}$] (Dinf) {};

\node[var-b,below left of=Dinf,label=right:HornSat] (Horn) {};
\node[var-b,below left of=Horn,label=right:] (K3) {};
\node[var-b,above left of=K3,label=right:{}] (3Lin2) {};
\node[above of=K3] (2SATn) {};
\node[above of=2SATn] (2SAT2) {};
\node[var-b,above of=2SAT2] (2SAT) {};

\node[var-b,above left of=2SAT,label=right:{}] (C2) {};

\path
    (C1) edge (P2)
    (P2) edge (B2)
    (B2) edge (B3)
    (B2) edge (D2)
    (D2) edge (D3)
    (B3) edge (D3)
    (Binf) edge (Dinf)
    (Dinf) edge (Horn)
    (Horn) edge (K3)
    (2SAT) edge (K3)
    (3Lin2) edge (K3)
    (C2) edge (3Lin2)
    (C2) edge (2SAT)
    (D2) edge (2SAT)
    (P2) edge (C2)
    ;

\end{tikzpicture}
}

\newcommand{\bA}{{\mathfrak A}}
\newcommand{\bB}{{\mathfrak B}}
\newcommand{\bC}{{\mathfrak C}}
\newcommand{\bD}{{\mathfrak D}}
\newcommand{\bF}{{\mathfrak F}}
\newcommand{\bG}{{\mathfrak G}}
\newcommand{\bI}{{\mathfrak I}}
\newcommand{\bK}{{\mathfrak K}}

\newcommand{\bP}{{\mathfrak P}}
\newcommand{\bT}{{\mathfrak T}}
\newcommand{\bZ}{{\mathfrak Z}}
\newcommand{\connects}{\ensuremath{\leftrightarrow^\ast}}
\DeclareMathOperator{\Csp}{CSP}
\DeclareMathOperator{\ar}{ar}

\title
{Symmetric Linear Arc Monadic Datalog \\
and Gadget Reductions}
\author{Manuel Bodirsky and Florian Starke\thanks{
Both authors have been funded by the European Research Council (Project POCOCOP, ERC Synergy Grant
101071674). Views and opinions expressed are however those of the authors only and do not necessarily reflect
those of the European Union or the European Research Council Executive Agency. Neither the European Union nor
the granting authority can be held responsible for them.
} \\ Institut f\"ur Algebra, TU Dresden}
\date{\today}

\begin{document}

\maketitle

\begin{abstract}
     A Datalog program \emph{solves} a constraint satisfaction problem (CSP) if and only if it derives the goal predicate precisely on the unsatisfiable instances of the CSP. 
    There are three Datalog fragments that are particularly important for finite-domain constraint satisfaction: \emph{arc monadic Datalog},
    \emph{linear Datalog}, and \emph{symmetric linear Datalog}, each having good computational properties. 
    We consider the fragment of Datalog where we impose all of these restrictions  simultaneously, i.e., we study \emph{symmetric linear arc monadic (slam) Datalog}. 
    We 
    characterise the CSPs that can be solved by a slam Datalog program as those that have a gadget reduction to a particular Boolean constraint satisfaction problem.
    We also present exact characterisations in terms of a homomorphism duality (which we call \emph{unfolded caterpillar duality}), and in universal-algebraic terms (using known minor conditions, namely the existence of \emph{quasi Maltsev operations} and \emph{$k$-absorptive operations of arity $nk$}, for all $n,k \geq 1$). 
    Our characterisations also imply that the question whether a given finite-domain CSP can be expressed by a slam Datalog program is decidable.
    
\end{abstract}


\tableofcontents

\section{Introduction}
   Datalog is an important concept linking database theory with the theory of constraint satisfaction. 
It is by far the most intensively studied formalism for polynomial-time tractability in constraint satisfaction. Datalog allows to formulate algorithms that are based on iterating local inferences, aka \emph{constraint propagation} or \emph{establishing local consistency}; this has been made explicit by Feder and Vardi in their groundbreaking work where they also formulate the finite-domain CSP dichotomy conjecture~\cite{FederVardi}. 
Following their convention, we say that a Datalog program $\Pi$ \emph{solves} a CSP if $\Pi$ derives the goal predicate on an instance of the CSP if and only if the instance is unsatisfiable.


The class of CSPs that can be solved by  Datalog is closed under so-called `gadget reductions' (a result due to Larose and Z\'adori~\cite{LaroseZadori}).
    In such a reduction, the variables in an instance of a constraint satisfaction problem are replaced by tuples of variables of some fixed finite length, and the constraints are replaced by \emph{gadgets} (implemented by conjunctive queries; a formal definition can be found in Section~\ref{sect:pp-constructions}); many of the well-known reductions between computational problems can be phrased as gadget reductions. Datalog is sufficiently powerful to simulate such gadget reductions; this has been formalised by Atserias, Bulatov, and Dawar in~\cite{AtseriasBulatovDawar} and the connection has been sharpened recently by Dalmau and Opr\v{s}al~\cite{DalmauOprsalLocal}. 

Feder and Vardi showed that Datalog cannot solve systems of linear equations over finite fields, even though such systems can be solved in polynomial time~\cite{FederVardi}. They suggest that the ability to simulate systems of linear equations should essentially be the only reason for a CSP to not be in Datalog. This conjecture was formalised by Larose and Z\'adori~\cite{LaroseZadori}: they observed that if systems of linear equations admit a gadget reduction to a CSP, then the CSP is not in Datalog, and they conjectured that otherwise the CSP can be solved by Datalog. 
This conjecture was proved by Barto and Kozik in 2009~\cite{BartoKozikFOCS09}, long before the resolution of the finite-domain CSP dichotomy conjecture by Bulatov~\cite{BulatovFVConjecture} and  by Zhuk~\cite{ZhukFVConjecture,Zhuk20}. 

 Datalog programs can be evaluated in polynomial time; but even a running time in $O(n^3)$ on a sequential computer can be prohibitively expensive in practise. This is one of the reasons why syntactic \emph{fragments} of Datalog  have been studied, which often come with better computational properties. 
 
 \subsection{Arc Monadic Datalog} 
    In \emph{monadic Datalog}, we restrict  the arity of the inferred predicates of the Datalog program  to one (i.e., all the \emph{IDBs} are monadic).
    In \emph{arc Datalog} we restrict each rule to a single input relation symbol (i.e., the body contains a single \emph{EDB};  for formal definitions, see Section~\ref{sect:datalog}). 
    
     An important Datalog fragment is \emph{arc monadic Datalog},
    which is still powerful enough to express the famous \emph{arc consistency procedure} in constraint satisfaction. 
     The arc consistency procedure has already been studied by Feder and Vardi~\cite{FederVardi}, and has many favorable properties: it can be evaluated in linear time and linear space. It is used as an important pre-processing step in the algorithms for both of the mentioned CSP dichotomy proofs, and it is also used in many practical implementations of algorithms in constraint satisfaction. The arc consistency procedure is still extremely powerful, and can for instance solve the P-complete HornSat Problem. 

Feder and Vardi characterised the power of the arc consistency procedure 
in terms of \emph{tree duality} (see Section~\ref{sect:dualities}), a natural combinatorial property which has been studied intensively in the graph homomorphism literature in the 90s (see, e.g.,~\cite{HNZ,HNBook}). Their characterisation has several remarkable consequences: one is that also the class of CSPs that can be solved by an arc monadic Datalog program is closed under gadget reductions. Another one is a \emph{collapse result} for Datalog when it comes to finite-domain CSPs, namely that monadic Datalog collapses to arc monadic Datalog: in fact, 
    if a finite-domain CSP can be solved by a monadic Datalog program, then it can already be solved by the arc consistency procedure (i.e., by a program in arc monadic Datalog). This statement is false without the restriction to finite-domain CSPs; in fact, there are infinite-domain CSPs that can be solved by a monadic Datalog program, but not by a program in arc monadic Datalog (Bodirsky and Dalmau~\cite{BodDalJournal}). 
 
\subsection{Linear Datalog} 
Besides arc monadic Datalog, there are other natural fragments of Datalog. The most notable one is \emph{linear Datalog}~\cite{Kanellakis}. 
Linear Datalog programs can be evaluated in non-deterministic logarithmic space (NL), and hence cannot express P-hard problems (unless P=NL). 
Dalmau~\cite{LinearDatalog} asked whether the converse is true as well, i.e., whether every finite-domain CSP which is in NL can be solved by a linear Datalog program~\cite{LinearDatalog}. 
This is widely treated as a conjecture, to which we refer as the \emph{linear Datalog conjecture}; it is one of the biggest open problems in finite-domain constraint satisfaction.

There are some sufficient 
conditions for solvability by linear Datalog (see Bulatov, Kozik, and Willard~\cite{BartoKozikWillard}
and Carvalho, Dalmau, and Krokhin~\cite{CarvalhoDalmauKrokhin})      and some necessary conditions (Larose and Tesson~\cite{LaroseTesson})
but the results {might} leave a gap.
    {A candidate for an example that falls into the gap can already be found among the CSPs of orientations of trees;} see Bodirsky, Bul\'in, Starke, and Wernthaler~\cite{otrees}.
Again, linear Datalog is closed under gadget reductions~\cite{StarkeDiss}.
And indeed, if HornSat has a gadget-reduction in a finite-domain CSP, then the finite-domain CSP cannot be solved by a linear Datalog program~\cite{AfratiCosmadakis}.

\subsection{Symmetric Linear Datalog} 
A further restriction is \emph{symmetric linear Datalog}, introduced by 
    Egri, Larose, and Tesson~\cite{EgriLaroseTessonLogspace}. 
    Symmetric linear Datalog programs can be evaluated in deterministic logspace~(L).
    Egri, Larose, and Tesson conjecture that every finite-domain CSP which is in L can be solved by a symmetric linear Datalog program~\cite{EgriLaroseTessonLogspace}; we refer to this conjecture as the \emph{symmetric linear Datalog conjecture}. 
    Symmetric linear Datalog is closed under gadget reductions~\cite{StarkeDiss}.
    Since directed reachability is not in symmetric linear Datalog~\cite{EgriLT08}, 
    it follows that every CSP that admits a gadget reduction from directed reachability cannot be solved by a symmetric linear Datalog program. 
    Egri, Larose, and Tesson also suggest that 
    this might be the only additional condition for containment in symmetric linear Datalog, besides the known necessary conditions to be in linear Datalog. 
    {A combinatorial characterisation of symmetric linear Datalog  has been presented by Egri~\cite{Egri14}.}
    
    Kazda~\cite{Kazda-n-permute} confirms the symmetric linear Datalog conjecture conditionally on the truth of the linear Datalog conjecture, i.e., he 
    shows that if a finite-domain CSP is in linear Datalog and does not admit a gadget reduction from a CSP that corresponds to the directed reachability problem, then 
    it is in symmetric linear Datalog (generalizing an earlier result of Dalmau and Larose~\cite{DalmauLarose08}). 
    

\subsection{Our Contributions} 
In this paper, we study the Datalog fragment that can be obtained by combining all the previously considered restrictions, namely \emph{symmetric linear arc monadic (slam) Datalog}. 
Before stating our result we illustrate this fragment with some examples. For $n \geq 1$, let $\bP_n$ be the directed path with $n$ vertices and $n-1$ edges. 
An example of a slam Datalog program which solves {the CSP of $\bP_2$} is 
    \begin{align*}
    A(x) & \;  {:}{-} \; 
    E(x,y) & 
    \text{goal} & \; {:}{-} \; E(x,y), A(y) 
    \end{align*}
     (in this case, the program is even recursion-free). 
    An example of a slam Datalog program which solves CSP$(\bP_3)$, this time with recursion and IDBs $A$ and $B$, is 
    \begin{align*}
        A(x) & \; {:}{-} \; E(x,y) &
        B(x) & \; {:}{-} \; A(y), E(x,y) \\
        A(y) & \; {:}{-} \; B(x), E(x,y)  & 
        \text{goal} & \; {:}{-} \;  B(y), E(x,y). 
    \end{align*}
The idea why this program is correct is that a finite digraph $\bA$ has a homomorphism to $\bP_3$ if and only if certain orientations of paths (those of \emph{net length three}; for a formal description, see Example~\ref{expl:zn}) do not have a homomorphism to $\bA$; and the program derives the goal predicate on $\bA$  precisely if there is a homomorphism from such a path to $\bA$.

    It follows from our results that the class of CSPs that can be solved by slam Datalog programs is closed under gadget reductions, despite the many restrictions that we imposed. 
We provide a \emph{full description} of the power of a Datalog fragment in terms of gadget reductions:  we show that a CSP can be solved by a slam Datalog program if and only if it has a gadget reduction to  {the CSP of $\bP_2$}.\footnote{The statement even holds for infinite-domain CSPs, since being solved by an arc monadic Datalog program implies the existence of a finite template~\cite{BodDalJournal} and admitting a gadget reduction to a finite-domain CSP implies the existence of a finite template as well~\cite{DalmauOprsalLocal}.}
The particular role of the structure $\bP_2$ is explained by the fact that
it is a representative of the unique class of CSPs which is non-trivial and \emph{weakest} with respect to gadget reductions -- a formalisation of this can be found in Section~\ref{sect:pp-constructions}.\footnote{We mention that $\bP_2$ is a Boolean structure which simultaneously satisfies the Schaefer conditions of being Horn, dual Horn, affine, and bijunctive~\cite{Schaefer}.} 
This shows that slam Datalog is the 
smallest non-trivial fragment of Datalog that is closed under gadget reductions. 

Our main result (Theorem~\ref{thm:main}) establishes a tight connection between the power of slam Datalog and 
various central themes in structural combinatorics and universal algebra. Specifically, the power of slam Datalog can be characterised by 
\begin{itemize}
    \item 
a new combinatorial duality which we call \emph{unfolded caterpillar duality} (restricting the concept of \emph{caterpillar duality} of Carvalho, Dalmau, and Krokhin~\cite{Lattice-Ops}), and using ideas that appear implicitly~\cite{DalmauLarose08,EgriLaroseTessonLogspace,Kazda-n-permute} {and explicitly~\cite{Egri14,DalmauLICS15}} in the literature on symmetric Datalog, and by 
\item 
the existence of a \emph{quasi Maltsev polymorphism} 
(a central concept in universal algebra) in combination with \emph{$kn$-ary $k$-absorbing polymorphisms for every $k,n \geq 1$} (introduced in ~\cite{Lattice-Ops} as well). 
\end{itemize}
Our result also implies that the following \emph{meta-problem} can be decided algorithmically: given a finite structure $\bB$, can {the CSP of $\bB$} be solved by a slam Datalog program? 

\medskip 
\subsection{Related Results}
Solvability of finite-domain CSPs by (unrestricted) Datalog 
was first studied by Feder and Vardi; they proved that {the CSP of $\bB$} can be solved by Datalog if and only if $\bB$ has \emph{bounded treewidth duality}, and they showed that CSPs for systems of linear equations over finite Abelian groups cannot be solved by Datalog. Larose and Zadori~\cite{LaroseZadori} showed that solvability by Datalog is preserved by gadget reductions and they asked whether having a gadget reduction from CSPs for systems of linear equations is not only a sufficient, but also a necessary condition for not being solvable by Datalog. This question was answered positive by Barto and Kozik~\cite{BoundedWidth}. 
Kozik, Krokhin, {Valeriote}, and Willard~\cite{Maltsev-Cond} gave a characterisation of Datalog in terms of minor conditions.

Linear (but not necessarily symmetric) monadic arc Datalog has been studied by Carvalho, Dalmau, and Krokhin~\cite{Lattice-Ops}; our proof builds on their result. 
In their survey on Datalog fragments and dualities in constraint satisfaction~\cite{BulatovKrokhinLarose08} 
Bulatov, Krokhin, and Larose 
write \emph{``it would be interesting to find (\dots) an appropriate notion of duality for symmetric (Linear) Datalog (...)''}. We do find such a notion for symmetric linear arc monadic  Datalog, namely unfolded caterpillar duality (Theorem~\ref{thm:main}). 

Another fragment of Datalog consists of the set of conjunctive queries; CSPs that can be solved by such Datalog programs are precisely the CSPs that are first-order expressible, by Rossman's theorem~\cite{Rossman08}. This is also known to be equivalent to the CSP having \emph{finite duality}. 
However, note that first-order definability is not preserved under gadget reductions (as we will see in Section~\ref{sect:dualities}).

\section{Preliminaries}
We write $[n]$ for the set $\{1,\dots,n\}$ and $[m,n]$ for the set $\{m,m+1,\dots,n\}$. 
We say that a tuple $a \in A^k$, for $k \in {\mathbb N}$, is \emph{injective} if $a$ is injective when viewed as a function from $[k]$ to $A$.  

\subsection{Structures and Graphs}
We assume familiarity with the concepts of relational structures and first-order formulas from mathematical logic, as introduced for instance in~\cite{Hodges}. 
The arity of a relation symbol $R$ is denoted by $\ar(R)$. If $\bA$ is a $\tau$-structure, 
then we sometimes use the same symbol for
$R \in \tau$ and the respective relation $R^\bA$ of $\bA$. 
We write $\bA[S]$ for the substructure of $\bA$ induced on $S$. 

A \emph{(directed) graph} is a relational structure with a single binary relation $E$. 
{If the edge relation is symmetric, then the graph is called \emph{undirected}.}  
For instance, the \emph{clique with $n$ vertices} is the {(undirected)} graph $\bK_n$ with domain $[n]$ and edges $E^{\bK_n} \coloneqq \{(a,b) \mid a \neq b\}$. 
Let $\bG$ be a graph. 
An \emph{(undirected) path from $a$ to $b$ in $\bG$} is a tuple $P = (a_1,\dots,a_n)$ such that $a_1,\dots,a_{n}$ are pairwise distinct, $a_1=a$, $a_n=b$, 
and for all $i \in [n-1]$ there is an edge between $a_i$ and $a_{i+1}$ (from $a_i$ to $a_{i+1}$ or from $a_{i+1}$ to $a_i$). 
If $i \in [2,n-1]$, then we say that $P$ \emph{passes through} $a_i$. 
An {undirected} graph $\bG$ is called \emph{connected} if for any two elements $a,b$ there exists a path from $a$ to $b$ in $\bG$.
A \emph{cycle} is a path from $a$ to $b$ of length $n$ at least three such that there is an edge between $a$ and $b$.
A graph is called \emph{acyclic} if it does not contain any cycle. 
A graph is called a \emph{tree} if it is connected and acyclic. 
A graph has \emph{girth $k$} if the length of the shortest cycle is $k$. 


{If $\bA$ is a $\tau$-structure and $E$ is an equivalence relation on the domain $A$ of $\bA$, then
$\bA/E$ denotes the $\tau$-structure $\bB$ whose domain is the set of equivalence classes $A/E$ of $E$, 
and where for every $R \in \tau$ of arity $n$, we have $(S_1,\dots,S_n) \in R^{\bB}$ if and only if $(a_1,\dots,a_n) \in R^{\bA}$ for some $a_1 \in S_1,\dots,a_n \in S_n$.} 

\subsection{Homomorphisms and CSPs}
Let $\tau$ be a relational signature and let $\bA$ and $\bB$ be $\tau$-structures. Then a \emph{homomorphism from $\bA$ to $\bB$} is a map $h \colon A \to B$ such that for $R \in \tau$, say of arity $k$, we have
$(h(a_1),\dots,h(a_k)) \in R^{\bB}$ whenever $(a_1,\dots,a_k) \in R^{\bA}$.
An embedding of $\bA$ into $\bB$ is an injective map
$e \colon A \to B$ such that
$(e(a_1),\dots,e(a_k)) \in R^{\bB}$ if and only if  $(a_1,\dots,a_k) \in R^{\bA}$.
We write $\bA \to \bB$ if there exists a homomorphism from $\bA$ to $\bB$ and $\bA \not\to \bB$ if there exists no homomorphism from $\bA$ to $\bB$.

If $\tau$ is a finite relational signature and $\bB$ is a $\tau$-structure, 
then $\Csp(\bB)$ denotes
the class of all finite $\tau$-structures $\bA$ such that $\bA \to \bB$. 
It can be viewed as a computational problem. For example, $\Csp(\bK_n)$ consists of the set of all finite $n$-colourable graphs, and can therefore be viewed as the $n$-colorability problem. 
Clearly, for finite structures $\bB$, this problem is always in NP. Note that from the database perspective, by the work of Chandra and Merlin~\cite{ChandraMerlin}, $\Csp(\bB)$ can be viewed as the {query} complexity of conjunctive queries over $\bB$.

A $\tau$-structure $\bB$ is \emph{homomorphically equivalent} to a $\tau$-structure $\bC$ if there are homomorphisms from $\bB$ to $\bC$ and vice versa. Clearly, homomorphically equivalent structures have the same CSP. 
A relational $\tau$-structure $\bC$ is called a \emph{core} 
if all endomorphisms of $\bC$ are embeddings. 
It is well-known and easy to see that every finite structure $\bB$ is homomorphically equivalent to a core $\bC$, 
and that all core structures $\bC$ that are homomorphically equivalent to $\bB$ are isomorphic; therefore, we refer to $\bC$ as \emph{the} core of $\bB$.


\subsection{Primitive Positive Constructions}
\label{sect:pp-constructions}
A $\tau$-formula $\phi$ is called a \emph{conjunctive query} (in constraint satisfaction and model theory such formulas are called \emph{primitive positive}, or short \emph{pp}) if it is built from atomic formulas (including atomic formulas of the form $x=y$) 
using only conjunction and existential quantification. 
If $\bB$ is a $\tau$-structure, and $\phi$ is a conjunctive query over the signature $\tau$, then 
$R := \{(t_1,\dots,t_k) \mid \bB \models \phi(t_1,\dots,t_k)\}$ is called the \emph{relation defined by $\phi$}. 

\begin{definition}\label{def:canonicalDatabase}
The \emph{canonical database} of a conjunctive query $\phi$ over the signature $\tau$ is the $\tau$-structure $\bB$ that can be constructed as follows: Let $\phi'$ be obtained from $\phi$ by renaming all existentially quantified variables such that no two quantified variables have the same name. 
Let $\phi''$ be obtained from $\phi'$ by removing all conjuncts of the form $x = y$ in $\phi'$ and by identifying variables $x$ and $y$ if there is a conjunct $x = y$ in $\phi'$.
Then $\bB$ is the $\tau$-structure whose domain is the set of variables of $\phi''$ such that for every $R \in \tau$ we have \[R^{\bB}=\{(v_1,\dots,v_k)\mid R(v_1,\dots,v_k)\text{ is a conjunct of }\phi''\}.\]
The \emph{canonical conjunctive query} of a structure $\bB$ with signature $\tau$ is the $\tau$-formula
with variables $B$ given by 
\[\bigwedge_{{R}\in\tau}\bigwedge_{t\in R^{\bB\hspace{-4mm}\phantom{(k)}}}R(t_1,\dots,t_{\ar(R)}).\]
\end{definition}

Observe that the canonical database of the canonical conjunctive query of a structure $\bB$ equals $\bB$. The following concepts have been introduced by Barto, Opr\v{s}al, and Pinsker~\cite{wonderland}. 

\begin{definition}
    A \emph{($d$-th) pp-power} of a $\tau$-structure $\bB$ is a structure $\bC$ with domain $B^d$ such that every relation of $\bC$ of arity $k$ is definable by a conjunctive query in $\bB$ as a relation of arity $dk$. 
    A structure has a \emph{primitive positive (pp) construction from $\bB$} if it is homomorphically equivalent to a pp-power of $\bB$. 
\end{definition}

Primitive positive constructions turned out to be \emph{the} essential tool for classifying the complexity of finite-domain CSPs, because if $\bC$ has a pp-construction from $\bB$, then there is a so-called \emph{gadget reduction} from $\Csp(\bC)$ to $\Csp(\bB)$; in fact, the converse is true as well, see Dalmau and Opr\v{s}al~\cite{DalmauOprsalLocal}. 

\begin{definition}\label{def:gadget-reduction}
Let $\mathcal B$ be a class of finite $\tau$-structures and let $\mathcal C$ be a class of finite $\rho$-structures. Then 
    a \emph{($d$-dimensional) gadget reduction} from ${\mathcal C}$ to ${\mathcal B}$ consists of  a conjunctive query $\phi_R$ of arity $dk$ over the signature $\tau$ for every $R \in \rho$ of arity $k$ { such that the following map $r$ from finite $\rho$-structures
    to finite $\tau$-structures satisfies $\bC\in\mathcal C$ if and only if $r(\bC)\in\mathcal B$ for all $\bC$.
    For a finite $\rho$-structures $\bC$ we obtain $r(\bC)$ as follows:}
    \begin{itemize}
        \item Replace each element $c$ of $\bC$ by the $d$-tuple $((c,1),\dots,(c,d))$.
        \item 
        {For every $R \in \rho$ of arity $k$ and every tuple $(t_1,\dots,t_k) \in R^{\bC}$, introduce a new element for every existentially quantified variable in $\phi_R$ and define relations for the relation symbols from $\tau$ such that the substructure induced by the new elements and $\{(t_1,1),\dots,(t_1,d),\dots,(t_k,1),\dots,(t_k,d)\}$ induce a copy of the canonical database of $\phi_R$ in the natural way, where $=$ is interpreted as a binary relation symbol (hence no identification of variables occurs). Let $\bC'$ be the resulting $\tau\cup\{=\}$-structure.}
        \item {Let $S$ be the smallest equivalence relation that contains $=^{\bC'}$. Then $r(\bC)$ is defined as the $\tau$-reduct of $\bC'/S$.}
    \end{itemize}
\end{definition}

For instance, the solution to the Feder-Vardi conjecture mentioned in the introduction states that $\Csp(\bB)$, for a finite structure $\bB$, is NP-hard if and only if $\bK_3$ has a pp-construction from $\bB$ (unless P=NP); 
by what we have stated above, 
this is true if and only if
the 3-coloring problem has a gadget reduction to $\Csp(\bB)$. 

\begin{example}\label{expl:pn}
    Let $\tau = \{E\}$ be the signature that consists of a single binary relation symbol $E$ whose elements we call \emph{edges}. 
    Let $\bP_n$ be the $\tau$-structure with the domain $\{1,2,\dots,n\}$ and edges $\{(1,2),(2,3),\dots,(n-1,n)\}$. 
    Then $\bP_2$ pp-constructs $\bP_3$, where the edge relation of $\bP_3$ is defined by \[\phi_E(x_{1,1},x_{1,2},x_{2,1},x_{2,2})\coloneqq E(x_{1,1},x_{2,2})\wedge (x_{1,2}=x_{2,1}).\] 
    Hence $\Csp(\bP_3)$ has a gadget reduction to $\Csp(\bP_2)$. Figure~\ref{fig:gadget-reduction-example} shows an example of the map $r$ from Definition~\ref{def:gadget-reduction}. 
\begin{figure}
    \centering
    \begin{tikzpicture}
        \node at (0,-2) {$\bC$};
        \node[scale=0.8] at (0,0) {\tikz{
        \node[var-b] (0a) at (0,0) {};
        \node[var-b] (1a) at (0,1) {};
        \node[var-b] (2a) at (0,2) {};
        \node[var-b] (1b) at (1,1) {};
        \node[var-b] (2b) at (1,2) {};
        \node[var-b] (3b) at (1,3) {};
        \path[->,>=stealth']
            (0a) edge (1a)
            (1a) edge (2a)
            (1b) edge (2a)
            (1b) edge (2b)
            (2b) edge (3b)
        ;
        }};

        \node at (1.5,0) {$\mapsto$};

        \node at (3,-2) {$\bC'$};
        \node[scale=0.8] at (3,0) {\tikz{
        \node[var-b] (0a0) at (0,0) {};
        \node[var-b] (1a0) at (0,1) {};
        \node[var-b] (2a0) at (0,2) {};
        \node[var-b] (1b0) at (1,1) {};
        \node[var-b] (2b0) at (1,2) {};
        \node[var-b] (3b0) at (1,3) {};
        \node[var-b] (0a1) at (0.4,0) {};
        \node[var-b] (1a1) at (0.4,1) {};
        \node[var-b] (2a1) at (0.4,2) {};
        \node[var-b] (1b1) at (1.4,1) {};
        \node[var-b] (2b1) at (1.4,2) {};
        \node[var-b] (3b1) at (1.4,3) {};
        \path[->,>=stealth']
            (0a0) edge (1a1)
            (1a0) edge (2a1)
            (1b0) edge (2a1)
            (1b0) edge (2b1)
            (2b0) edge (3b1)
        ;
        \path[dashed]
            (0a1) edge (1a0)
            (1a1) edge (2a0)
            (1b1) edge (2a0)
            (1b1) edge (2b0)
            (2b1) edge (3b0)
        ;
        }};
        
        \node at (4.5,0) {$\mapsto$};

        \node at (6,-2) {$r(\bC)$};
        \node[scale=0.7] at (6,0) {\tikz{
        \node[var-b] (0a0) at (0,0) {};
        \node[var-b] (1a0) at (0,1) {};
        \node[var-b] (1b0) at (1,1) {};
        \node[var-b] (2b0) at (1,2) {};
        \node[var-b] (2a1) at (0.4,2) {};
        \node[var-b] (2b1) at (1.4,2) {};
        \node[var-b] (3b1) at (1.4,3) {};
        \path[->,>=stealth']
            (0a0) edge (2b0)
            (1a0) edge (2a1)
            (1b0) edge (2a1)
            (1b0) edge (2b1)
            (2b0) edge (3b1)
        ;}};
    \end{tikzpicture}
    \caption{{An instance $\bC$ of $\Csp(\bP_3)$ that is reduced to $r(\bC)$, which is an instance of $\Csp(\bP_2)$. The dashed edges in the diagram of $\bC'$ represent $=^{\bC'}$.}} 
    \label{fig:gadget-reduction-example}
\end{figure}
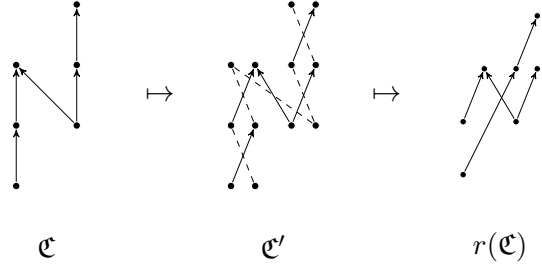
\end{example}

\begin{remark}\label{rem:p2}
    It is known that pp-constructibility is transitive~\cite{wonderland}, and
    the corresponding poset has a largest element (which is represented by $\bP_1$), and
    all other elements are below the element which is represented by $\bP_2$; see, e.g., {Propositions~3.2.2 and~3.2.5 in~\cite{VucajDiss}}.
\end{remark}

Given the fundamental importance of conjunctive queries and homomorphisms in database theory, we believe that pp-constructions and gadget reductions are an interesting concept for database theory as well.

\subsection{Datalog}
\label{sect:datalog}
Let $\tau$ and $\rho$ be finite relational signatures such that $\tau \subseteq \rho$. 
A \emph{Datalog program} is a finite set of rules of the form
$$ \phi_0 \; {:}{-} \; \phi_1,\dots,\phi_n$$
where each $\phi_i$ is an atomic $\tau$-formula. The formula $\phi_0$ is called the \emph{head} of the rule, and the sequence $\phi_1,\dots,\phi_n$ is called the \emph{body} of the rule. The symbols in $\tau$ are called \emph{EDBs} (\emph{extensional database predicates}) and the other symbols from $\rho$ are called \emph{IDBs} (\emph{intensional database predicates}). In the rule heads, only IDBs are allowed. There is one special IDB of arity 0, which is called the \emph{goal predicate}. 
IDBs might also appear in the rule bodies. We view the set of rules as a recursive specification of the IDBs in terms of the EDBs -- for a detailed introduction, see, e.g.,~\cite{BodDalJournal}. 
A Datalog program is called 
\begin{itemize}
    \item \emph{linear} if in each rule, at most one IDB appears in the body (we then assume without loss of generality that in every rule whose body contains an IDB, the IDB is listed first). 
\item \emph{arc} if each rule involves at most one EDB. 
\item \emph{symmetric} if it is linear and for every rule $\phi_0 \; {:}{-} \; \phi_1,\phi_2,\dots,\phi_n$ where $\phi_0$ and $\phi_1$ are {atomic formulas whose relation symbol is an IDB}, 
the Datalog program also contains the \emph{reversed rule}  $\phi_1 \; {:}{-} \; \phi_0,\phi_2,\dots,\phi_n$.\footnote{Note that we do not have to exclude that $\phi_0$ is the goal predicate, because 
we may always add its symmetric version without changing the set of structures on which the goal predicate is derived.} 
\end{itemize}

If $\bB$ is a $\tau$-structure, then 
we say that $\Csp(\bB)$ is \emph{solved} by a Datalog program $\Pi$ with EDBs $\tau$ if the following holds: the goal predicate is derived by $\Pi$ on a finite $\tau$-structure $\bA$ if and only if there is \emph{no} homomorphism from $\bA$ to $\bB$. 

We say that a Datalog program has \emph{width $(\ell,k)$}
if all IDBs have arity at most $\ell$, and if every rule has at most $k$ variables. 
For given $(\ell,k)$ and a structure $\bB$, there exists a Datalog program $\Pi$ of width $(\ell,k)$ with the remarkable 
property that if some Datalog program of width $(\ell,k)$ solves $\Csp(\bB)$, then $\Pi$ solves $\Csp(\bB)$. 
This Datalog program is referred to as the \emph{canonical Datalog program for $\bB$ of width $(\ell,k)$}, and is constructed as follows~\cite{FederVardi}:
For every relation $R$ over $B$ of arity at most $\ell$, we introduce a new IDB. The empty relation of arity 0 plays the role of the goal predicate. 
Then $\Pi$ contains all rules 
$\phi \; {:}{-} \; \phi_1,\dots,\phi_n$
with at most $k$ variables 
such that the formula 
$\forall \bar x (\phi_1 \wedge \dots \wedge \phi_n \Rightarrow \phi)$ holds in the expansion of $\bB$ by all IDBs. 
If the canonical Datalog program  for $\bB$  derives the goal predicate on a finite structure $\bA$, then 
there is no homomorphism from $\bA$ to $\bB$ (see, e.g.,~\cite{BodDalJournal}). 

If $k$ is the maximal arity of the EDBs, 
we may restrict the canonical Datalog program of width $(1,k)$ to those rules with only unary IDBs and at most one EDB; in this case, we obtain the canonical arc monadic Datalog program, 
 which is also known as the \emph{arc consistency procedure}. 
Analogously, we may define the canonical \emph{linear}, 
and the canonical \emph{symmetric} Datalog program. We may also combine these restrictions, and in particular obtain 
a definition the \emph{canonical slam Datalog program}, i.e., the canonical symmetric linear arc monadic Datalog program, which has not yet been studied in the literature before. 

The following lemma can be shown analogously to the well-known fact for unrestricted canonical Datalog programs of width $(\ell,k)$ (see, e.g.,~\cite{BodDalJournal}). 

\begin{lemma}\label{lem:can-sound}
Let $\bB$ be a finite structure with a finite relational signature,
and let $\Pi$ be the canonical slam Datalog program for $\bB$. 
If $\bA$ is a finite structure with a homomorphism to $\bB$, then 
$\Pi$ does not derive the goal predicate on $\bA$. 
\end{lemma} 

\subsection{The Incidence Graph} 
Several results from graph theory concerning 
acyclicity and high girth can be generalised to general structures. To formulate these generalisations, we need the concept of an  
 \emph{incidence graph} 
 of a relational structure $\bA$. 
 
 
 The \emph{incidence graph}
of a structure $\bA$ with the relational signature $\tau$ is the bipartite {undirected} 
graph where one colour class is $A$, and the other consists of all pairs of the form $(t,R)$ such that $t \in R^{\bA}$ and $R \in \tau$ {(in some texts, the incidence graph is a \emph{multigraph}; we prefer not to introduce multigraphs, but will comment on the difference later when it becomes relevant)}. 
We put an edge between $a$ and $(t,R)$ if $t_i=a$ for some $i$. 
The \emph{girth} of an (undirected) graph $\bG$ is the length of the shortest cycle in $\bG$. 
We say that a relational structure is a \emph{generalised tree}
if its incidence graph is a tree. 
A \emph{leaf} of a generalised tree $\bT$ is an element of $T$ which has degree one in the incidence graph of $\bT$.


A structure $\bB$ is called \emph{injective} if all tuples that are in some relation in $\bB$ are injective (i.e., have no repeated  entries; {for such structures, the incidence multigraph used in some texts has no multiple edges, and hence can be viewed as a usual undirected graph}). 
A structure is called an \emph{(injective) tree} if it is injective and its incidence graph is a tree; {note that this is standard in the literature~\cite{LLT}.}

\begin{theorem}[Sparse incomparability lemma for structures~\cite{FederVardi}]
\label{thm:sparse}
Let $\tau$ be a finite relational signature. 
    Let $\bA$ and $\bB$ be $\tau$-structure with finite domains such that $\bA\not\to\bB$. Then for every $m \in {\mathbb N}$ there exists an injective finite structure $\bA'$ whose incidence graph has girth at least $m$, such that  $\bA'\to\bA$ and $\bA'\not\to\bB$. 

\end{theorem}


\subsection{Dualities}
\label{sect:dualities}
For a $\tau$-structure $\bB$ and a class of $\tau$-structures $\mathcal F$ 
the pair $({\mathcal F},\bB)$ is called a \emph{duality pair} if a finite structure $\bA$ has a homomorphism to $\bB$ if and only if no structure $\bF \in {\mathcal F}$ has a homomorphism to $\bA$. 
Several forms of duality pairs will be relevant here, depending on the class of structures ${\mathcal F}$.

A $\tau$-structure $\bB$ has \emph{finite duality} if there exists a finite set of $\tau$-structures ${\mathcal F}$ such that $({\mathcal F},\bB)$ is a duality pair. The property of having finite duality is among the very few notions studied in the context of constraint satisfaction which is \emph{not} preserved under gadget reductions, as illustrated in the following example.

\begin{example}\label{expl:zn}
    As in Example~\ref{expl:pn}, let $\tau = \{E\}$.
    {For $n\geq 1$, let $\bZ_{n-1}$ be the $\tau$-structure with domain 
    $\{s,1,1',2,2',\dots,n,n',t\}$ and edges 
    \[\{(s,1),(1,1'),(2,1'),(2,2'),\dots,(n,(n-1)'),(n,n'),(n',t)\}.\]}
    Then 
    \begin{itemize}
        \item $\bP_2$ has finite duality, witnessed by the duality pair $(\{\bP_3\},\bP_2)$, 
    \item recall from Example~\ref{expl:pn} that $\bP_2$ pp-constructs $\bP_3$~\cite{maximal-digraphs}, so $\Csp(\bP_3)$ has a gadget reduction to $\Csp(\bP_2)$, 
    but 
    \item $\bP_3$ does not have finite duality: this is witnessed by the fact that $(\{\bZ_n \mid n \in {\mathbb N}\},\bP_3)$ is a duality pair~\cite{HNBook}, and that there is no homomorphism from $\bZ_n$ to $\bZ_m$ for $n < m$.  \qedhere 
    \end{itemize}
\end{example}

\begin{example}\label{expl:b2}
    Let $\rho = \{E,Z\}$ be the signature that consists of a binary relation symbol $E$ and a unary relation symbol $Z$. Let $\bB_2$ be the structure with domain $\{0,1\}$ where
    \begin{align*}
        E^{\bB_2} & \coloneqq \{(1,1),(0,1),(1,0)\} \\
        Z^{\bB_2} & \coloneqq \{0\} 
    \end{align*}
    Let $\bP_2'$ be the $\rho$-expansion of $\bP_2$ where $Z^{\bP_2'} \coloneqq \{0,1\}$. 
    Then $(\{\bP_2'\},\bB_2)$ is a duality pair. 
\end{example}

A more robust form of duality is \emph{tree duality}, which plays a central role in constraint satisfaction, and is studied in the graph homomorphism literature in the 90s. A structure $\bB$ has \emph{tree duality} 
if there exists a (not necessarily finite) set of trees ${\mathcal F}$ such that $({\mathcal F},\bB)$ is a duality pair. 
The following is well known; see Theorem 7.4 in~\cite{PCSP}. The equivalence of 1.\ and 3.\ is from~\cite{FederVardi}; also see~\cite{GraphHomomorphisms}. 

\begin{theorem}\label{thm:ac1}
    Let $\bB$ be a finite $\tau$-structure. 
    Then the following are equivalent:
    \begin{enumerate}
        \item $\bB$ has tree duality; 
        \item $\bB$ has a pp-construction from
        $(\{0,1\};\{0\},\{1\},\{0,1\}^3 \setminus \{(1,1,0)\})$; 
        \item $\bB$ can be solved by arc consistency. 
    \end{enumerate}
\end{theorem}

There are finite structures with tree duality that have a P-complete CSP, such as the structure
$(\{0,1\};\{0\},\{1\},\{0,1\}^3 \setminus \{(1,1,0)\})$, which is essentially the Boolean HornSAT problem. In the following, we therefore introduce more restrictive forms of dualities. 
\begin{figure}
    \centering
    \begin{tikzpicture}[rotate=90,yscale=1.3]
        \node[var-b,label=right:1] (1) at (0,0) {};
        \node[var-b,label=right:{$((1,2,3),R)$}] (R123) at (1,0) {};
        \node[var-b,label=below right:3] (3) at (2,0) {};
        \node[var-b,label=below:2] (2) at (1,1) {};

        \node[var-b,label=right:{$((3,6),E)$}] (R36) at (3,0) {};
        \node[var-b,label=right:6] (6) at (4,0) {};
        
        \node[var-b,label=below:{$((3,4,5),R)$}] (R345) at (2,1) {};
        \node[var-b,label=below:4] (4) at (2,2) {};
        \node[var-b,label=right:5] (5) at (3,1) {};

        \node[var-b,label=below:{$((3),P)$}] (R3) at (2,-1) {};
        \node[var-b,label=right:{$((1),P)$}] (R1) at (-1,0) {};
        \node[var-b,label=right:{$((6),P)$}] (R6) at (5,0) {};
        
        \path 
            (1) edge (R123)
            (2) edge (R123)
            (3) edge (R123)
            (3) edge (R345)
            (4) edge (R345)
            (5) edge (R345)
            (3) edge (R36)
            (6) edge (R36)
            (3) edge (R3)
            (1) edge (R1)
            (6) edge (R6);
    \end{tikzpicture}\hspace{5mm}
    \begin{tikzpicture}[rotate=90,yscale=1.1]
        \node[var-b,label=right:1] (1) at (0,0) {};
        \node[var-b,label=right:{$((1,2),E)$}] (R12) at (1,0) {};
        \node[var-b,label=right:2] (2) at (2,0) {};

        \node[var-b,label=right:{$((2,3,4),R)$}] (R234) at (3,0) {};
        \node[var-b,label=right:4] (4) at (4,0) {};

        \node[var-b,label=below:3] (3) at (3,1) {};
        \node[var-b,label=below:{$((3),P)$}] (R3) at (3,2) {};

        \node[var-b,label=right:{$((4,5),E)$}] (R45) at (5,0) {};
        \node[var-b,label=right:5] (5) at (6,0) {};
        
        \path 
            (1) edge (R12)
            (2) edge (R12)
            (2) edge (R234)
            (3) edge (R234)
            (4) edge (R234)
            (4) edge (R45)
            (5) edge (R45)
            (3) edge (R3);
    \end{tikzpicture}
    \caption{An example of the incidence graph of a caterpillar (left) and of a structure that is \textbf{not} a caterpillar (right; {the trouble is that $((3),P)$ is not adjacent to the vertical path, and that this path is the only candidate for the spine}).}
    \label{fig:exampleCaterpillar}
\end{figure}
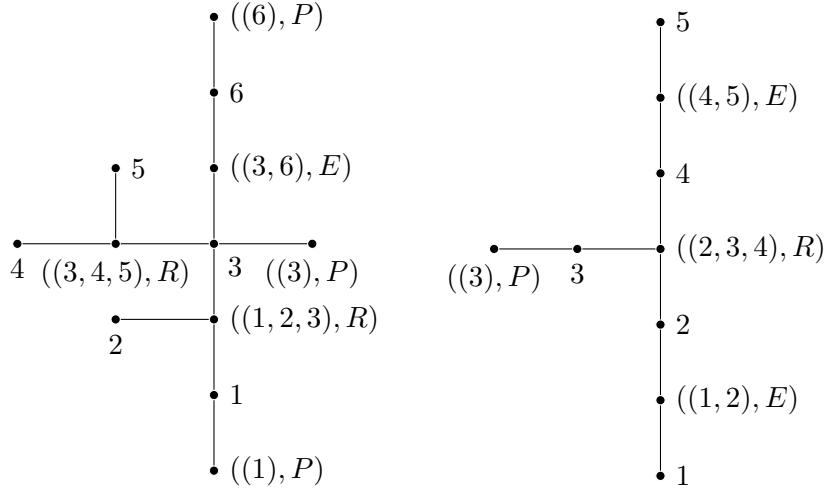

\begin{definition}
A relational structure $\bA$ is called a \emph{generalised caterpillar} if 
\begin{itemize}
    \item it is a generalised tree,
\item its incidence graph $\bG$ contains a path $P= (a_1,\dots,a_n)$, {called the \emph{spine} of the caterpillar},  such that every vertex in $G \setminus \{a_1,\dots,a_n\}$ of the form $(t,R)$ is adjacent (in $\bG$) to a vertex in $P$. 
\end{itemize}
\end{definition}
See Figure~\ref{fig:exampleCaterpillar} (left) for an example. 
    A relational structure $\bA$ is called an \emph{(injective) caterpillar} if it is injective and a generalised caterpillar
    (this definition of a caterpillar is equivalent to the one  given in~\cite{Lattice-Ops}). 
    A structure $\bB$ has \emph{caterpillar duality} if there exists 
a set of caterpillars ${\mathcal F}$ such that $({\mathcal F},\bB)$ is a {duality pair~\cite{CarvalhoDK08}}.
The structures $\bB$ with caterpillar duality have been characterised by~\cite{Lattice-Ops} (Theorem~\ref{thm:caterpillar}) in terms of linear arc Datalog. 

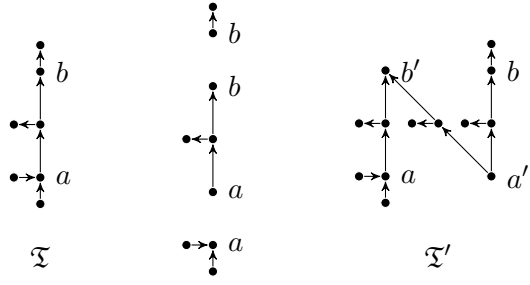
\begin{figure}
    \centering
    \begin{tikzpicture}[scale=0.7,rotate=90]
    
        \node at (-.5,0) {$\bT$};
        \node[var-b] (0) at (0.5,0) {};
        \node[var-b] (1) at (1,0.5) {};
        \node[var-b, label={right:$a$}] (a) at (1,0) {};
        \node[var-b] (2) at (2,0) {};
        \node[var-b] (3) at (2,0.5) {};
        \node[var-b, label={right:$b$}] (b) at (3,0) {};
        \node[var-b] (4) at (3.5,0) {};

        \path[->,>=stealth']
            (0) edge (a)
            (1) edge (a)
            (a) edge (2)
            (2) edge (3)
            (2) edge (b)
            (b) edge (4)
            ;
    \end{tikzpicture}\hspace{12mm}
        \begin{tikzpicture}[scale=0.7,rotate=90]
        \node[var-b] (0) at (0.5,0) {};
        \node[var-b] (1) at (1,0.5) {};
        \node[var-b, label={right:$a$}] (a) at (1,0) {};

        \node[var-b, label={right:$a$}] (a') at (2,0) {};
        \node[var-b] (2) at (3,0) {};
        \node[var-b] (3) at (3,0.5) {};
        \node[var-b, label={right:$b$}] (b) at (4,0) {};

        \node[var-b, label={right:$b$}] (b') at (5,0) {};
        \node[var-b] (4) at (5.5,0) {};

        \path[->,>=stealth']
            (0) edge (a)
            (1) edge (a)
            (a') edge (2)
            (2) edge (3)
            (2) edge (b)
            (b') edge (4)
            ;
    \end{tikzpicture}\hspace{12mm}
    \begin{tikzpicture}[scale=0.7,rotate=90]
            \node at (-.5,-1) {$\bT'$};
        \node[var-b] (0) at (0.5,0) {};
        \node[var-b] (1) at (1,0.5) {};
        \node[var-b, label={right:$a$}] (a) at (1,0) {};
        \node[var-b] (2) at (2,0) {};
        \node[var-b] (3) at (2,0.5) {};
        \node[var-b, label={right:$b'$}] (b) at (3,0) 
        {};
        
        \node[var-b, label={right:$a'$}] (a') at (1,-2) {};
        \node[var-b] (2') at (2,-1) {};
        \node[var-b] (3') at (2,-0.5) {};
        
        \node[var-b] (2'') at (2,-2) {};
        \node[var-b] (3'') at (2,-1.5) {};
        \node[var-b, label={right:$b$}] (b'') at (3,-2) 
        {};
        \node[var-b] (4) at (3.5,-2) {};

        \path[->,>=stealth']
            (0) edge (a)
            (1) edge (a)
            (a) edge (2)
            (2) edge (3)
            (2) edge (b)
            (2') edge (b)
            (2') edge (3')
            (a') edge (2')
            (a') edge (2'')
            (2'') edge (3'')
            (2'') edge (b'')
            (b'') edge (4)
            ;
    \end{tikzpicture}

    \caption{Example of an $(a,b)$-unfolding $\bT'$ of a tree $\bT$.}
    \label{fig:unfoldingExample}
\end{figure}

To capture the power of \emph{symmetric} linear arc Datalog, we present a more restrictive form of duality. 
Let $\bT$ be an (injective) 
tree and 
let $a$ and $b$ be distinct elements of $\bT$.
Write the canonical query of $\bT$ as $\phi_a \wedge \phi_{a,b} \wedge \phi_b$ where
$\phi_a$ contains all conjuncts of the form 
$R(\bar u)$ such that in the incidence graph, 
{the unique path from the vertex $(\bar u,R)$ to the vertex $b$ passes through $a$.} Similarly, we define $\phi_b$, switching the roles of $a$ and $b$. 
Note that $\phi_a$ and $\phi_b$ do not share any conjuncts. 
All the remaining conjuncts of the canonical query form {$\phi_{a,b}$}. 
Let {$\psi_a$ be obtained from $\phi_a$ by existentially quantifying all variables except for $a$.
Similarly, we define $\psi_b$, switching the roles of $a$ and $b$.}
Let $\psi$ be obtained from $\phi_{a,b}$ by existentially quantifying all variables except for $a$ and $b$. 
%
{The following concept is closely related to notions
introduced by Egri~\cite{Egri14} and used in~\cite{DalmauLICS15} (we discuss the differences below);} see Figure~\ref{fig:unfoldingExample} for an example.

\begin{definition}\label{def:unfold}
Let $\bT$ be an (injective) tree and let $a$ and $b$ be distinct elements of $\bT$ that are not leaves. 
The \emph{$(a,b)$-unfolding of $\bT$} is 
the canonical database of the formula
$${\psi_a(a) \wedge \psi(a,b') \wedge \psi(a',b') \wedge \psi(a',b) \wedge \psi_b(b)}$$
{where $\psi_a$, $\psi_b$, and $\psi$ obtained from the canonical database of $\bT$ as explained above.} 
An \emph{unfolding} of $\bT$ is a structure $\bT'$ that is obtained by a sequence $\bT =\bT_1,\bT_2,\dots,\bT_n = \bT'$, {where for all $i \in \{1,\dots,n-1\}$ there are $a_i,b_i$ in $\bT_{i}$  such that
$\bT_{i+1}$ is an $(a_i,b_i)$-unfolding of $\bT_{i}$.}
\end{definition}

Note that an unfolding of a tree $\bT$ is again a tree and has a homomorphism to $\bT$. 
It can also be shown that the unfolding of a caterpillar is a caterpillar as well (Lemma~\ref{lem:caterpillar-unfolded}). 
We say that a structure $\bB$ has \emph{unfolded caterpillar duality} if there exists a set $\mathcal F$ of caterpillars such that $({\mathcal F},\bB)$ is a duality pair, and
${\mathcal F}$ contains every unfolding of a caterpillar in ${\mathcal F}$. 
Clearly, unfolded caterpillar duality implies caterpillar duality. 

{The difference to $(j,k)$-symmetric bounded pathwidth duality introduced by Egri~\cite{Egri14} and to \emph{zigzag-realizations} from~\cite{DalmauLICS15} is that in contrast to those works we aim for a characterisation of \emph{arc} Datalog programs, which means that our obstruction sets need to be caterpillars rather than $(j,k)$-paths, which makes our results formally incomparable.} 


\subsection{Minor Conditions}
If $\bB$ is a structure and $k \geq 1$, then a \emph{polymorphism of $\bB$ of arity $k$} is a homomorphism from $\bB^k$ to $\bB$. The set of all polymorphisms of $\bB$ is denoted by $\Pol(\bB)$. 
An \emph{(operation) clone} is a set of operations which contains the projections and is closed under composition. 
Note that $\Pol(\bB)$ is a clone. 
An operation $f \colon B^n \to B$ is called \emph{idempotent} if $f(x,\dots,x) = x$ for all $x \in B$. A clone is called \emph{idempotent} if all of its operations are idempotent.

If $\bA$ and $\bB$  are a structures and $k \geq 1$, then a \emph{polymorphism of $(\bA,\bB)$ of arity $k$} is a homomorphism from $\bA^k$ to $\bB$. The set of all polymorphisms of $(\bA,\bB)$ is denoted by $\Pol(\bA,\bB)$. Let $f\colon A^n\to B$ be a function and let $\sigma\colon [n]\to [m]$, then the map
\begin{align*}
    f_\sigma\colon A^m&\to B\\
    (a_1,\dots,a_m)&\mapsto f(a_{\sigma(1)},\dots,a_{\sigma(n)})
\end{align*}
is called a \emph{minor} of $f$. 
A \emph{minion} is a set or functions from $A^n$ to $B$ that is closed under taking minors.
Note that $\Pol(\bA,\bB)$ is a minion. 
Let $\mathscr M$ and $\mathscr N$ be minions. A map $\xi\colon \mathscr M \to \mathscr N$ is called a \emph{minion homomorphism} if $\xi$ preserves arity and for every $f\in\mathscr M$ of arity $n$ and every map $\sigma\colon [n]\to [m]$ we have 
\[\xi(f_\sigma)=(\xi(f))_\sigma.\]

Let $\tau$ be a \emph{function signature}, i.e., a set of function symbols, each equipped with an arity. 
A \emph{minor condition} is a finite set $\Sigma$ of 
{\emph{height-one identities}}, i.e., expressions of the form 
\begin{equation*}
    f(x_1,\dots,x_n) \approx g(y_1,\dots,y_m)
\end{equation*}
where $f$ is an $n$-ary function symbol from $\tau$,
$g$ is an $m$-ary function symbol from $\tau$, 
and $x_1,\dots,x_n,y_1,\dots,y_m$ are (not necessarily distinct) variables. 
If $\mathscr M$ is a minion, then a map $\xi \colon \tau \to \mathscr M$ \emph{satisfies} a minor condition $\Sigma$ if 
for every {height-one identity} $f(x_1,\dots,x_n) \approx g(y_1,\dots,y_m) \in \Sigma$ and for every assignment $s \colon \{x_1,\dots,x_n,y_1,\dots,y_n\} \to B$ we have $$\xi(f)(s_1(x_1),\dots,s(x_n)) = \xi(g)(s(y_1),\dots,s(y_m)).$$ 
We say that a minion $\mathscr M$
\emph{satisfies} $\Sigma$ if 
there exists a map $\xi \colon \tau \to \mathscr M$ 
that satisfies $\Sigma$.\footnote{It is convenient and standard practise to notationally drop the distinction between $f \in \tau$ and $\xi(f) \in {\mathscr M}$.}
If $\Sigma$ and $\Sigma'$ are minor conditions, then we say that {
\emph{$\Sigma$ implies $\Sigma'$ (for clones)}} if every clone that satisfies $\Sigma$ also satisfies $\Sigma'$. {Note that in general $\Sigma$ implies $\Sigma'$ for clones does not mean that every minion which satisfies  $\Sigma$ also satisfies $\Sigma'$.} 
We present some concrete minor conditions that are relevant in the following. 

\begin{definition}
    An operation $m \colon B^3 \to B$ is called a \emph{quasi Maltsev operation} if it satisfies
    the minor condition 
     $$m(x,x,y) \approx m(y,x,x) \approx m(y,y,y).$$
    A \emph{Maltsev operation} is an idempotent quasi Maltsev operation. A \emph{quasi minority operation} is a quasi Maltsev operation $m$ that additionally satisfies $$m(x,y,x) \approx {m(y,y,y)}$$ and a \emph{minority operation} is an idempotent quasi minority operation.
\end{definition}

    \begin{definition}
    An operation $f \colon B^m \to B$ is called 
    a \emph{quasi majority operation}
     if it satisfies the minor condition 
     $$m(x,x,y) \approx m(x,y,x) \approx m(y,x,x) \approx m(x,x,x).$$
    A \emph{majority operation} is an idempotent quasi majority operation. 
\end{definition}

The following was shown in in~\cite{CCC,FederVardi} {(also see Theorem~7.1 in~\cite{OligoClone}).}

\begin{proposition}\label{prop:decomp}
    Let $\bB$ be a finite relational $\tau$-structure {which is a core}. Then the following are equivalent:
    \begin{itemize}
        \item $\bB$ has a quasi majority polymorphism. 
        \item Every relation with a primitive positive definition in $\bB$ has a definition by a conjunction of primitive positive formulas, each with at most two free variables.  
    \end{itemize}
\end{proposition}

\begin{example}\label{expl:Bn}
    Generalising Example~\ref{expl:b2}, the structure $\bB_n$ has domain $\{0,1\}$, signature $\{{\boldsymbol 0},R_n\}$ and 
    the relations $\{0\}$ and 
    $\{0,1\}^n \setminus \{(0,\dots,0)\}$. Define the structure $\bF_n$ with domain $\{1,\dots,n\}$, signature $\{{\boldsymbol 0},R_n\}$, ${\boldsymbol 0}\coloneqq\{1,\dots,n\}$, and $R_n\coloneqq \{(1,\dots,n)\}$. Note that $(\{\bF_n\},\bB_n)$ is a duality pair and that $\bF_n$ is a tree, but not a caterpillar.
    The structure $\bB_n$ has finite duality, but no quasi Maltsev polymorphism and no 
    quasi majority polymorphism. It follows from Theorems~\ref{thm:caterpillar} and \ref{thm:main} that $\Csp(\bB)$ can be solved by linear arc monadic Datalog but not by slam Datalog. 
\end{example}

\begin{definition}
    Let $k,n \in {\mathbb N}_{>0}$. 
    An operation $f \colon B^{kn} \to B$ is called \emph{$k$-block symmetric} if it satisfies the following condition 
\begin{align}
f(x_{11},\dots,x_{1k},\dots,x_{n1},\dots,x_{nk}) \approx f(y_{11},\dots,y_{1k},\dots,y_{n1},\dots,y_{nk}) \label{eq:blocksym}
\end{align}
whenever $\{S_1,\dots,S_n\} = \{T_1,\dots,T_n\}$ where $S_i = \{x_{i1},\dots,x_{ik}\}$ and $T_i = \{y_{i1},\dots,y_{ik}\}$.
If $k=1$ or $n=1$ then $f$ is called \emph{totally symmetric}.

If $f$ is $k$-block symmetric and
$S_1,\dots,S_n$ are subsets of $B$ of size at most $k$, 
then we also write $f(S_1,\dots,S_n)$ instead of $f(x_{11},\dots,x_{1k},\dots,x_{n1},\dots,x_{nk})$ where $\{x_{i1},\dots,x_{ik}\} = S_i$. We say that $f$ is \emph{$k$-absorptive} if it satisfies 
$$f(S_1,S_2,\dots,S_n) \approx f(S_2,S_2,S_3,\dots,S_n)$$ whenever $S_2 \subseteq S_1$. 
\end{definition}

\begin{remark}\label{rem:majo}
    Note that every structure with a $2$-absorptive polymorphism $f$ of arity $6$ also has the quasi majority polymorphism $m$ given by 
    $$m(x,y,z) \coloneqq f(x,y,z,x,y,z),$$ 
    because 
    $\{x\} \subseteq \{x,z\}$ and hence 
    $$m(x,x,z) {{}\approx{}} f(x,x,z,x,x,z) {{}\approx{}} f(x,x,x,x,x,x) {{}\approx{}} m(x,x,x)$$ and similarly $$m(x,z,x) {{}\approx{}} m(z,x,x) {{}\approx{}} m(x,x,x).$$ 
\end{remark}

The list of equivalent statements from Theorem~\ref{thm:ac1} can now be extended as follows.  

\begin{theorem}[\cite{DalmauPearson,FederVardi}]
\label{thm:ac2}
    Let $\bB$ be a finite $\tau$-structure. 
    Then 
        $\bB$ has tree duality if and only if 
    $\bB$ has totally symmetric polymorphisms of all arities. 
\end{theorem}

We will make crucial use of the following theorem. 

\begin{theorem}[Theorem 16 in \cite{Lattice-Ops}]\label{thm:caterpillar}
    Let $\bB$ be a finite relational $\tau$-structure. 
    Then the following are equivalent. 
    \begin{enumerate}
        \item \label{item:cp-1} $\bB$ has caterpillar duality. 
        \item $\Csp(\bB)$ can be solved by a linear arc monadic Datalog program. 
        \item $\Pol(\bB)$ contains for every $k,n \geq 1$ an $k$-absorbing operation of arity  $kn$. 
    \item \label{item:cp-4}
    $\bB$ is homomorphically equivalent to a structure $\bB'$ with binary polymorphisms $\sqcup$ and $\sqcap$ such that $(B',\sqcup,\sqcap)$ is a (distributive) lattice. 
    \end{enumerate}
\end{theorem}


\subsection{Indicator Structures}
\label{sect:indicator}
In this section we revisit a common theme in constraint satisfaction, the concept of an \emph{indicator structure} of a minor condition. 
To simplify the presentation, we only define the indicator structure for minor conditions with only one 
function symbol. For our purposes, this is without loss of generality, because for clones over a finite domain, every minor condition is equivalent to such a restricted minor condition. If $f \colon C^n \to C$ and $g \colon C^m \to C$ are operations, then the \emph{star product} $f * g$ is defined to be the operation defined as
$$ (x_{1,1},\dots,x_{n,m}) \mapsto f(g(x_{1,1},\dots,x_{1,m}),\dots,g(x_{n,1},\dots,x_{n,m})).$$

\begin{lemma}\label{lem:one-f}
    Let $\Sigma$ be a minor condition. Then there exists a minor condition $\Sigma'$ with a single function symbol such that 
    a clone over a finite domain satisfies $\Sigma$ if and only if it satisfies $\Sigma'$. 
\end{lemma}
\begin{proof}
    First note that 
    for every clone ${\mathscr D}$ on a finite set there exists an \emph{idempotent} clone ${\mathscr C}$ on a finite set 
    which is equivalent to it with respect to minion homomorphisms,
    i.e., there are minion homomorphism
    from ${\mathscr D}$ to ${\mathscr C}$ and vice versa. {(If all unary operations in ${\mathscr D}$ are bijective, then we may take for ${\mathscr C}$ the set of all idempotent operations in ${\mathscr D}$.)} 
    It is well-known and easy to see that if $f_1,\dots,f_n$ are the function symbols {of arities $k_1,\dots,k_n$, respectively,} that appear in $\Sigma$, 
    and $\mathscr C$ satisfies $\Sigma$, then
    $\mathscr C$ also 
    contains an operation $g$ of arity $m$
    such that for every $i \in [n]$ there exists {$\alpha_i \colon [m] \to [k_i]$} such that 
    $g_{\alpha_i} = f_i$ (use that ${\mathscr C}$ is closed under the star product and idempotent). 
    Note that ${\mathscr C}$ satisfies a 
    {height-one identity}
    $(f_i)_{\beta} \approx (f_j)_{\gamma}$
    if and only if
    ${\mathscr C}$ satisfies a 
    {height-one identity}
    $(g_{\alpha_i})_{\beta} \approx (g_{\alpha_j})_{\gamma}$. 
\end{proof}


\begin{definition}\label{def:indicatorStructure}
Let $\bB$ be a relational $\tau$-structure and let $\Sigma$ be a minor condition with a single function symbol $f$ of arity $m$. 
Let $\sim$ be the smallest equivalence relation on $B^m$ such that 
 $a \sim b$ if $\Sigma$ contains $f(x_1,\dots,x_m)\approx f(y_1,\dots,y_m)$
such that there is a map $s \colon \{x_1,\dots,x_m,y_1,\dots,y_m\} \to B$ with 
$a = (s(x_1),\dots,s(x_m))$ and $b = (s(y_1),\dots,s(y_m))$. {Note that if $a\sim b$, then $\Sigma$ implies $f(a)\approx f(b)$, here the entries of $a$ and $b$ are viewed as variables.}
Then the \emph{indicator structure of $\Sigma$ with respect to $\bB$} is the $\tau$-structure $\bB^m/_{\sim}$. 
\end{definition}

The following is straightforward from the definitions. 

\begin{lemma}\label{lem:ind} 
Let $\bB$ be a structure and $\Sigma$ be a minor condition with a single function symbol $f$. Then 
$\bB$ has a polymorphism satisfying $\Sigma$ if and only if the indicator structure of $\Sigma$ with respect to $\bB$ has a homomorphism to $\bB$.
\end{lemma}

\section{Results}
In this section we state and prove our main result (Theorem~\ref{thm:main}), which characterises the power of slam Datalog in many different ways, including descriptions in terms of pp-constructability in $\bP_2$, minor conditions, 
unfolded caterpillar duality, and homomorphic equivalence to a structure with both lattice and quasi Maltsev polymorphisms.


\begin{theorem}\label{thm:main}    
Let $\bB$ be a structure with a finite domain and a finite relational signature $\tau$. 
Then the following are equivalent. 
    \begin{enumerate}
        \item \label{maltsev}
        $\Pol(\bB)$ contains a quasi Maltsev operation and $k$-absorptive operations
        of arity $nk$, for all $n,k \geq 1$. 
        \item 
        \label{can-sym-lin-arc}
        The canonical slam  Datalog program for $\bB$ solves $\Csp(\bB)$.
        \item \label{sym-lin-arc}
        Some slam Datalog program solves $\Csp(\bB)$.
        \item 
        \label{caterpillar}
        $\bB$ has unfolded caterpillar duality. 
        \item \label{non-degenerate-minor}
        If $\Pol(\bB)$ does not satisfy a minor condition $\Sigma$,
        then $\Sigma$ implies $f(x) \approx f(y)$. 
        \item \label{minor-p2} Every minor condition that holds in $\Pol(\bP_2)$ also holds in $\Pol(\bB)$. 
        \item \label{minion-hom} There is a minion homomorphism from $\Pol(\bP_2)$ to $\Pol(\bB)$. 
        \item \label{pp-p2} There is a pp-construction of 
        $\bB$ in $\bP_2$.
        \item \label{lattice} $\bB$ is homomorphically equivalent to a structure $\bB'$ such that $\Pol(\bB')$ contains a quasi Maltsev operation and operations $\sqcup$ and $\sqcap$ such that $(B',\sqcup,\sqcap)$ forms a (distributive) lattice.
    \end{enumerate}
    Moreover, if one of these items holds, {and $\bB$ is a core,} then there exists a structure $\bB'$ with a binary relational signature such that  $\Pol(\bB') = \Pol(\bB)$,
    and all the statements hold for $\bB'$ in place of $\bB$ as well. 
\end{theorem}

We first prove the equivalence of \eqref{maltsev}-\eqref{minor-p2} in cyclic order. We then explain how the equivalence of 
\eqref{minor-p2}-\eqref{pp-p2} follows from 
general results in the literature, and finally show the equivalence of~\eqref{maltsev} and~\eqref{lattice}. 
The proof of the theorem stretches over the following subsections.

\begin{example}
\label{expl:Tn}
    The structure $\bT_n$ is the transitive tournament with $n$ vertices, i.e., it has the 
    domain $[n]$ and the binary relation $<$. 
    Note that $\bT_2$ equals $\bP_2$. 
    It is easy to see that $(\{\bP_{n+1}\},\bT_n)$ is a duality pair. Since $\bP_{n+1}$ is a caterpillar Theorem~\ref{thm:caterpillar} implies that $\Csp(\bT_n)$ can be solved by a linear arc monadic Datalog program. However, 
    $\bT_n$ does not have a quasi Maltsev polymorphism for $n \geq 3$,
    and hence Theorem~\ref{thm:main} implies that $\Csp(\bT_n)$ cannot be solved by slam Datalog. 
\end{example}

{
It is well known {and easy to see} that a finite structure $\bB$ is pp-interconstructable with $\bP_1$, i.e., there is a pp-construction of $\bB$ in $\bP_1$ and vice versa, if and only if the core of $\bB$ has one element. Hence, by Remark~\ref{rem:p2}, we can use Theorem~\ref{thm:main} to characterises all finite structures that are pp-interconstructable with $\bP_2$. 
\begin{corollary}\label{cor:characterisationP2}
    Let $\bB$ be a structure with a finite domain and a finite relational signature. Then $\bB$ is pp-interconstructable with $\bP_2$ if and only if the core of $\bB$ has at least two elements and $\bB$ satisfies one of the items from Theorem~\ref{thm:main}.
\end{corollary}
}

\subsection{Symmetrizing Linear Arc Monadic Datalog}
The following lemma is used for the implication  from (\ref{maltsev}) to (\ref{can-sym-lin-arc}) in the proof of Theorem~\ref{thm:main}.
Note that in the canonical linear arc monadic Datalog program $\Pi$ we can use the `strongest possible rules'\footnote{These comments are intended to illustrate the challenges in the proof of next lemma; it will not be necessary to formalise what we mean by strongest possible rules.} when deriving the goal predicate. However, the canonical slam Datalog program $\Pi_S$ might need to use weaker rules in order to be able to apply symmetric rules later on in the derivation. See Example~\ref{exa:SLAMmustUseWeakRules}. 

\begin{lemma}\label{lem:MaltimpliesLAMDatalogIsSLAMDatalog}
Let $\bB$ be a finite structure with relational signature $\tau$ such that $\Pol(\bB)$ contains a quasi Maltsev operation. 
Let $\Pi$ be the canonical linear arc monadic Datalog program for $\bB$ and $\Pi_S$ be the canonical slam Datalog program for $\bB$. Then $\Pi$ can derive the goal predicate on a finite 
$\tau$-structure 
$\bA$ if and only if $\Pi_S$ can derive the goal predicate on $\bA$.
\end{lemma}
\begin{proof}
    Note that every rule of $\Pi_S$ is also a rule of $\Pi$. Hence if $\Pi_S$ can derive the goal predicate on $\bA$, then so can $\Pi$.
    {Consider a derivation of the goal predicate $G$ for the Datalog program $\Pi$ on $\bA$. Let $R_0,\dots, R_{n+1}$ be the rules of $\Pi$ used in this derivation. Since $\Pi$ is linear, we can use the suggestive notation 
    \[\vdash_{R_0} P_0(a_0)\vdash_{R_1} \cdots \vdash_{R_{n}}P_n(a_n)\vdash_{R_{n+1}}G\]
    for this derivation.}
    For $i\in[n]$ define the primitive positive formula $\Phi_i$ as follows. 
    Suppose that the rule $R_i$ is of the form
    $P_{i}(y) \;{:}{-}\; \Psi_i(x_1,\dots,x_k)\wedge P_{i-1}(x)$ for some atomic formula $\Psi_i$.
    We may assume that both $x$ and $y$ are among the variables $x_1,\dots,x_k$; otherwise, the canonical database of $P_{i}(y)\wedge\Psi_i(x_1,\dots,x_k)\wedge P_{i-1}(x)$ 
    is not connected. 
    Since $\Pi$ solves a CSP we may assume without loss of generality that such rules $R_i$ were not used in the derivation of the goal predicate.
    
    \begin{itemize}
        \item If $x\neq y$, then $\Phi_i(x,y)$ is obtained from $\Psi_{i}(x_1,\dots,x_k)$ by existentially quantifying all variables except for $x$ and $y$,
        \item otherwise, $x=y$ and $\Phi_i(x,x')$  is obtained from $\Psi_i(x_1,\dots,x_k)$ by existentially quantifying all variables except for $x$ and adding the conjunct $x=x'$.
    \end{itemize}
    Define $\Phi_0(x)$ and $\Phi_{n+1}(x)$ from the rules $R_0$ and $R_{n+1}$ in a similar fashion.
    For $i\in[n]$ define the binary relation $\to_i$ on $B$ that contains all tuples $(b,b')$ such that $\bB\models \Phi_i(b,b')$.
    We may assume without loss of generality that \[P_0=\{b\in B\mid \bB\models \Phi_0(b)\} =\Phi_0^{\bB}\] and that for all $b'\in P_{i+1}$ there exists a $b\in P_i$ such that $b\to_i b'$. In particular, this implies that $P_0$ is pp-definable.

    Let $Q_0,\dots,Q_n\subseteq B$ be the smallest sets such that
    \begin{itemize}
        \item $P_i\subseteq Q_i$,
        \item $b\in Q_i$ and $b\to_i b'$ implies $b'\in Q_{i+1}$, and
        \item $b'\in Q_{i+1}$ and $b\to_i b'$ implies $b\in Q_{i}$.
    \end{itemize}
    Note that there can be $b\in Q_i$ such that there is no $b'\in B$ with $b\to_i b'$. Analogously, there can be $b'\in Q_{i+1}$  such that there is no $b\in B$ with $b\to_i b'$.
    Let $\tilde R_0,\dots,\tilde R_{n+1}$ be the rules obtained from $R_0,\dots,R_{n+1}$ by replacing each occurrence of $P_i$ by $Q_i$ for all $i\in[0,n]$. 
    We will now show that 
    \[\vdash_{\tilde R_0} Q_0(a_0)\vdash_{\tilde R_1}\dots\vdash_{\tilde R_{n}}Q_n(a_n)\vdash_{\tilde R_{n+1}}G\] is a derivation of $\Pi_S$ on $\bA$.
    It suffices to show that $\tilde R_i$ is a rule of $\Pi_S$ for all $i\in[0,n+1]$.
    By definition $P_0\subseteq Q_0$. Hence $\bB\models\forall x(\Phi_0(x)\Rightarrow P_0(x))$ implies $\bB\models\forall x(\Phi_0(x)\Rightarrow Q_0(x))$. Therefore, $\tilde R_0$ is a rule of $\Pi_S$.

    Let $i\in[n]$. To show that {$\bB\models\forall x,y \big ((\Phi_i(x,y)\wedge Q_{i-1}(x))\Rightarrow Q_{i}(y) \big)$,} let $b\in Q_{i-1}$ and $b'\in B$ be such that $\bB\models \Phi_{i}(b,b')$. Then $b\to_{i} b'$ and $b'\in Q_{i}$ by the definition of $Q_{i}$. Analogously, we show that  {$\bB\models\forall x,y \big ((\Phi_{i}(x,y)\wedge Q_{i}(y))\Rightarrow Q_{i-1}(x) \big)$}. Therefore, $\tilde R_i$ is a rule of $\Pi_S$. 

To prove that $\tilde R_{n+1}$ is a rule in $\Pi_S$ we will show that for all $b_n\in Q_n$ we have $\bB\not\models\Phi_{n+1}(b_n)$. Let $b_n\in Q_n$ and assume that $\bB\models\Phi_{n+1}(b_n)$. 
    By the definition of $Q_0,\dots,Q_n$ and the condition that for all $b'\in P_{i+1}$ there exists a $b\in P_i$ such that $b\to_i b'$ we know that there exists a $b_0\in P_0$ such that $b_0$ and $b_n$ are connected by the symmetric transitive closure of $\bigcup_{i=1}^{n}\to_i$.  
    
    Let $m$ be a quasi Maltsev polymorphism of $\bB$. Applying   $m$ repeatedly to the connection of $b_0$ and $b_{n}$ as indicated in the following picture:
    
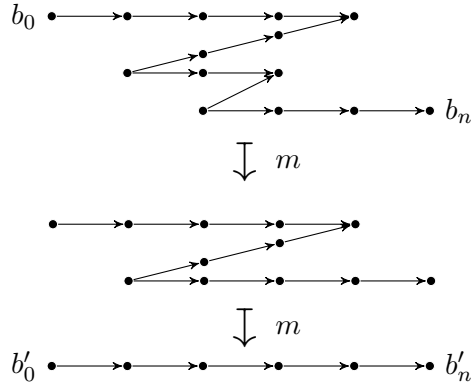
\begin{center}
        
    \begin{tikzpicture}
    
    \node at (0,0) {\tikz[yscale=0.25]{
        \node[var-b,label=left:$b_0$] at (0,0) (0) {};
        \node[var-b] at (1,0) (1) {};
        \node[var-b] at (2,0) (2) {};
        \node[var-b] at (3,0) (3) {};
        \node[var-b] at (4,0) (4) {};
        
        \node[var-b] at (1,-3) (1') {};
        \node[var-b] at (2,-2) (2') {};
        \node[var-b] at (3,-1) (3') {};

        \node[var-b] at (2,-3) (2'') {};
        \node[var-b] at (3,-3) (3'') {};

        \node[var-b] at (2,-5) (2''') {};
        \node[var-b] at (3,-5) (3''') {};
        \node[var-b] at (4,-5) (4''') {};
        \node[var-b,label=right:$b_n$] at (5,-5) (5''') {};
        
        \path[->,>=stealth']
            (0) edge (1)
            (1) edge (2)
            (2) edge (3)
            (3) edge (4)
            (3') edge (4)
            (2') edge (3')
            (1') edge (2')
            (1') edge (2'')
            (2'') edge (3'')
            (2''') edge (3'')
            (2''') edge (3''')
            (3''') edge (4''')
            (4''') edge (5''')
            ;
            }};
        \node[rotate=-90,label=above:$m$, scale=1.5] at (0,-1.3) {$\mapsto$};
    \node at (0,-2.5) {\tikz[yscale=0.25]{
        \node[var-b] at (0,0) (0) {};
        \node[var-b] at (1,0) (1) {};
        \node[var-b] at (2,0) (2) {};
        \node[var-b] at (3,0) (3) {};
        \node[var-b] at (4,0) (4) {};
        
        \node[var-b] at (1,-3) (1') {};
        \node[var-b] at (2,-2) (2') {};
        \node[var-b] at (3,-1) (3') {};

        \node[var-b] at (2,-3) (2'') {};

        \node[var-b] at (3,-3) (3''') {};
        \node[var-b] at (4,-3) (4''') {};
        \node[var-b] at (5,-3) (5''') {};
        
        \path[->,>=stealth']
            (0) edge (1)
            (1) edge (2)
            (2) edge (3)
            (3) edge (4)
            (3') edge (4)
            (2') edge (3')
            (1') edge (2')
            (1') edge (2'')
            (2'') edge (3''')
            (3''') edge (4''')
            (4''') edge (5''')
            ;
            }};
        \node[rotate=-90,label=above:$m$, scale=1.5] at (0,-3.5) {$\mapsto$};
    \node at (0,-4) {\tikz[yscale=0.25]{
        \node[var-b,label=left:$b'_0$] at (0,0) (0) {};
        \node[var-b] at (1,0) (1) {};
        \node[var-b] at (2,0) (2) {};
        \node[var-b] at (3,0) (3) {};
        \node[var-b] at (4,0) (4) {};
        \node[var-b,label=right:$b'_n$] at (5,0) (5''') {};
        
        \path[->,>=stealth']
            (0) edge (1)
            (1) edge (2)
            (2) edge (3)
            (3) edge (4)
            (4) edge (5''')
            ;
            }};
    \end{tikzpicture}
        \end{center}

    we can conclude that there exist $b'_0,b'_1,\dots,b'_n\in B$ such that \[b'_0\to_1 b'_1\to_2\dots\to_n b'_n\] and $\bB\models\Phi_{n+1}(b'_n)$. Since $P_0$ is pp-definable, it is preserved by $m$ and therefore $b'_0\in P_0$. 
    Hence, $b'_2\in P_2,\dots,b'_n\in P_n$. This contradicts that $R_{n+1}$ is a rule of $\Pi$ (as $\bB\models \Phi_{n+1}(b'_n)\wedge P_n(b'_n)$). Hence, $\bB\not\models\Phi_{n+1}(b_n)$ and $\tilde R_{n+1}$ is a rule of $\Pi_S$ as desired. 
\end{proof}

\begin{remark}
Observe that the canonical database of \[\Phi_0(x_0)\wedge\Phi_1(x_0,x_1)\wedge\dots\wedge\Phi_n(x_{n-1},x_{n})\wedge\Phi_{n+1}(x_{n})\] is a generalised caterpillar. 
\end{remark}

\begin{example}\label{exa:SLAMmustUseWeakRules}
    Consider the structure $\bB$ with domain $\{0,0',1,a,b,b'\}$, binary relation $E=\{(0,1),(0',1),(a,b),(a,b')\}$, and {constants for $0$ and $b$}. Let $\bA$ be the structure with domain $\{0,b\}$, the binary relation $E=\{(0,b)\}$, and all constants. Clearly, $\bA\not\to\bB$. The  canonical linear arc monadic Datalog program $\Pi$ for $\bB$ can derive the goal predicate using the 
    derivation {(in the same notation as in the proof of Lemma~\ref{lem:MaltimpliesLAMDatalogIsSLAMDatalog})}\[\vdash_{R_0}\{0\}(0)\vdash_{R_1}\{1\}(b)\vdash_{R_2} G.\] The  canonical slam Datalog program $\Pi_S$ for $\bB$ can also derive the goal predicate on $\bA$. It cannot use the rule $R_1$, because $R_1$ is not symmetric. However, it may use a different rule $\tilde R_1$:
    \[\vdash_{\tilde R_0}\{0,0'\}(0)\vdash_{\tilde R_1}\{1\}(b)\vdash_{R_2} G.\]
    Note that $R_0$ is also a rule of $\Pi_S$ but in order to apply $\tilde R_1$ the program needs to use the rule $\tilde R_0$ which is weaker than the rule $R_0$ (in the sense that the derived IDB is a strict superset). 
\end{example}

\subsection{Proving Unfolded Caterpillar Duality}

This section is devoted to the proof 
of the implication 
\eqref{sym-lin-arc} to~\eqref{caterpillar} in Theorem~\ref{thm:main}. 
We first prove a general result 
about obstruction sets for finite-domain CSPs that closes a gap in the presentation of the proof of Lemma~21 in~\cite{Lattice-Ops} and essentially follows from the sparse incomparability lemma (Theorem~\ref{thm:sparse}); we thank V\'ictor Dalmau for clarification. 

\begin{lemma}\label{lem:inj}
    Let $\bB$ be a finite structure and let $\mathcal F$ be a class of finite structures such that $(\mathcal F,\bB)$ is a duality pair. Define 
    \[\mathcal F'\coloneqq \{\bF\in\mathcal F \mid \text{$\bF$ is injective}\}.\]
    Then $(\mathcal F',\bB)$ is a duality pair.
\end{lemma}
\begin{proof}
Let $\tau$ be the signature of $\bB$; let $m$ be the maximal arity of the relation symbols in $\tau$. 
    Let $\bA$ be a finite $\tau$-structure 
    which does not homomorphically map to $\bB$. 
    By Theorem~\ref{thm:sparse}, there exists an injective finite structure $\bA'$ whose incidence graph has girth at least
    $3$, and which homomorphically maps to $\bA$ but not to $\bB$. 
    There exists a $\bF \in {\mathcal F}$ which homomorphically maps to $\bA'$. 
    Since $\bA'$ is injective, so is $\bF$. 
    It follows that $\bF \in {\mathcal F'}$. This implies that $({\mathcal F}',\bB)$ is a duality pair. 
\end{proof}

Note that Theorem~\ref{thm:caterpillar} implies that if $\Csp(\bB)$ is solved by  a linear arc monadic Datalog program, then $\bB$ has caterpillar duality; the proof given in~\cite{Lattice-Ops} only shows 
\emph{generalised caterpillar duality}, {i.e., there exists a  set $\mathcal F$ of generalised caterpillars such that $(\mathcal F,\bB)$ is a duality pair. However, Lemma~\ref{lem:inj} implies that
these notions are equivalent; similarly, the definition of tree duality and the analogously defined notion of generalised tree duality are equivalent.}

In order to prove the implication 
\eqref{sym-lin-arc} to~\eqref{caterpillar} in Theorem~\ref{thm:main}, it only remains to prove that $\bB$ also has \emph{unfolded} caterpillar duality (Lemma~\ref{lem:SLAMimpliesUnfoldedCaterpillar}). 
We first prove the following lemma, which has already been mentioned in Section~\ref{sect:dualities}.

\begin{lemma}\label{lem:caterpillar-unfolded}
    An unfolding of a caterpillar is a caterpillar as well.
\end{lemma}
\begin{proof}
    Let $\bD$ be a caterpillar and let $\bD'$ be an $(a,b)$-unfolding of $\bD$ for two non-leafs $a,b \in D$. Let $P = (a_1,\dots,a_n)$ be a longest possible path in the incidence graph of $\bD$ which shows that $\bD$ is a caterpillar. Note that since $P$ is longest possible it must pass through all non-leafs of $\bD$ (see Figure~\ref{fig:exampleCaterpillar}). In particular, $P$ passes  through $a$ and through $b$;
    without loss of generality it can be written as $$(\bar u,a,\bar v,b,\bar w).$$ 
    Let $\phi_1$, $\phi_2$, and $\psi$ be obtained from $\bD$ as in the definition 
    of the $(a,b)$-unfolding of $\bD$, so that 
    $\bD'$ is the canonical databases of $\phi_1(a) \wedge \psi(a,b') \wedge \psi(a',b') \wedge \psi(a',b) \wedge \phi_2(b)$. 
    Note that 
    \begin{itemize}
        \item 
        $(\bar u,a)$ is a path in the incidence graph of the canonical database of $\phi_1(a)$,
        \item
    $(a,\bar v,b)$
    is a path in the incidence graph of the canonical database of $\psi(a,b)$, 
    and 
    \item $(b,\bar w)$ is a path in the incidence graph of the  canonical database of $\phi_2(b)$.
    \end{itemize}
    Furthermore, each of these paths witnesses that the corresponding canonical database is a caterpillar.
    Let $(a,\bar v_1,b')$, 
    $(a',\bar v_2,b')$, 
$(a',\bar v_3,b)$
be the corresponding paths in 
the incidence graph of the canonical database of 
$\psi(a,b')$, $\psi(a',b')$, and $\psi(b',b)$, respectively. Let $\tilde v_2$ be $\bar v_2$ in reversed order. 
Then 
$$P' \coloneqq (\bar u,a,\bar v_1,b',\tilde v_2,a',\bar v_3,b,\bar w)$$ is a path in 
 the incidence graph $\bG'$ of $\bD'$. 
We claim that $P'$ witnesses that $\bD'$ is a caterpillar. 
This follows from the observation that if  
$\bC_i$, for $i \in \{1,2\}$, is a caterpillar with witnessing path $P_i=(\bar u_i,a_i)$ such that $a_i \in C_i$, then the structure obtained by taking the disjoint union of $\bC_1$ and $\bC_2$ and identifying $a_1$ and $a_2$ is a caterpillar,  witnessed by the path $(\bar u_1,a_1,\bar u_2)$. 

    The statement for unfoldings in general follows by induction. 
\end{proof}

\begin{lemma}\label{lem:SLAMimpliesUnfoldedCaterpillar}
    If $\Csp(\bB)$ is solved by a slam Datalog program $\Pi$, then $\bB$ has unfolded caterpillar duality. 
\end{lemma}
\begin{proof}
    As we have explained above, 
    Theorem~\ref{thm:caterpillar} (in combination with Lemma~\ref{lem:inj}) implies that there exists a set of caterpillars $\mathcal F$ such that 
    $(\mathcal F,\bB)$ is a duality pair. 
    Let ${\mathcal F}'$ 
    be the closure of ${\mathcal F}$ by all unfoldings of caterpillars in ${\mathcal F}$. 
    By Lemma~\ref{lem:caterpillar-unfolded}, we have that ${\mathcal F}'$ is a set of caterpillars as well. It remains to show that $({\mathcal F}',\bB)$  is a duality pair. Since $\mathcal F\subseteq\mathcal F'$ it suffices to show that no element in $\mathcal F'$ maps homomorphically to $\bB$. 
Let $\bF\in\mathcal F$, $a,b \in F$, and $\bF'$ be the $(a,b)$-unfolding of $\bF$ (Definition~\ref{def:unfold}). Let $\phi_1$, $\phi_2$, and $\psi$ be obtained from $\bF$ as in the definition of the $(a,b)$-unfolding of $\bF$ so that
    $\bF'$ is the canonical database of $\phi_1(a) \wedge \psi(a,b') \wedge \psi(a',b') \wedge \psi(a',b) \wedge \phi_2(b)$. Assume without loss of generality that for every existentially quantified variable $v$ in $\psi(a,b')$ the corresponding variables in $\psi(a',b')$ and $\psi(b',a)$ are $v'$ and $v''$, respectively. 
    Let $P_0(v_0)\vdash_{R_1}\dots \vdash_{R_n}P_n(v_n)$ be a derivation of $\Pi$ on $\bF$ (assuming that $\Pi$ already derived $P_0$ on $a$) with $v_0=a$, $v_n=b$, and $v_1,\dots,v_{n-1}\in F\setminus\{a,b\}$ such that $P_{i-1}$ occurs in the body of $R_i$ for every $i\in[n]$. Then $v_1,\dots,v_n$ are in the canonical database of $\psi(a,b)$. Hence $v_i$, $v'_i$, and $v''_i$ are elements in $\bF'$ for every $i\in[n]$. For every rule $R_i$ let $\tilde R_i$ denote the reversed rule, which is also a rule of $\Pi$ since $\Pi$ is symmetric. Note that 
    \begin{align*}
  P_0(a)\vdash_{R_1} \cdots \vdash_{R_n}P_n(b') & \vdash_{\tilde R_n} P_{n-1}(v'_{n-1}) \\ \cdots & \vdash_{\tilde R_1}P_0(a')\vdash_{R_1}P_1(v''_1) 
  \cdots  \vdash_{R_n}P_n(b)
    \end{align*}
    is a derivation of $\Pi$ on $\bF'$. 
    Therefore,  any derivation $d$ of $\Pi$ on $\bF$ deriving the IDB $P$ on an element $v$ can be transformed into a derivation $d'$ of $\Pi$ on $\bF'$ such that 
    \begin{itemize}
        \item  $d'$ derives the IDB $P$ on $v$ if $v$ is in the canonical databases of $\phi_1(a)$ or $\phi_2(b)$ and 
        \item  $d'$ derives the IDB $P$ on $v$ or on $v''$ if $v$ is in the canonical database of $\psi(a,b)$ and $v\notin\{a,b\}$. 
    \end{itemize}
    Since $\bF\in\mathcal F$ and $\Pi$ solves $\Csp(\bB)$, we have that $\Pi$ derives the goal predicate on $\bF$. Therefore, $\Pi$ derives the  goal predicate on $\bF'$ as well. Hence,  $\bF'$ does not map homomorphically to $\bB$.   
    It follows by induction that no unfolding of an element in $\mathcal F$ maps homomorphically to $\bB$. Therefore, $({\mathcal F}',\bB)$  is a duality pair.
\end{proof}

\subsection{Using Unfolded Caterpillar Duality}
This section proves the implication \eqref{caterpillar} to~\eqref{non-degenerate-minor} in Theorem~\ref{thm:main}.

\begin{lemma}\label{lem:interesting}
    Let $\bB$ be a relational structure with   unfolded caterpillar duality. If $\Pol(\bB)$ does not satisfy a minor condition $\Sigma$, then $\Sigma$ implies $f(x) \approx f(y)$. 
\end{lemma}
\begin{proof} 
Let $\tau$ be the signature of $\bB$. 
Suppose that 
{$\Pol(\bB)\not\models\Sigma$}. 
By Lemma~\ref{lem:one-f}, we may assume that $\Sigma$ only involves a single function symbol $f$ of arity $K$. 
Let $\bI$ be the indicator structure of $\Sigma$ with respect to $\bB$. Then $\bI\not\to\bB$ (Lemma~\ref{lem:ind}). 
By the caterpillar duality of $\bB$, there must be a caterpillar $\bC$ with a homomorphism to $\bI$ that does not have a homomorphism to $\bB$. 
In the following we will show that either there exists an unfolding of $\bC$ that has a homomorphism to $\bB$, which contradicts the unfolded caterpillar duality of $\bB$, or that $\Sigma$ implies $f(x)\approx f(y)$, which concludes the proof.


Let $P$ be a path witnessing that $\bC$ is a caterpillar. We may assume that $P$ is of the form 
\[(v_0,(s_1,R_1),\dots,(s_n,R_n),v_n).\]
This assumption is without loss of generality: 
if $P$ instead starts as follows 
$$((s_1,R_1),v_1,\dots )$$ 
then either $R_1$ is not unary, and we may prepend to $P$ an entry $v_0$ of $s_1$ which is different from $v_1$.  
Or $R_1$ is unary, in which case $(s_1,R_1)$ is a leaf (in the incidence graph of $\bC$) and we may simply discard this first vertex of $P$. In both cases, the new path still witnesses that $\bC$ is a caterpillar.
Analogously, we can ensure that $P$ does not end in a vertex of the form $(s_n,R_n)$. 

{Let $m$ be the maximal number of neighbours of any $v_i$ in the incidence graph of $\bC$ that are not in the path $P$. }
Define the primitive positive formulas $\phi_1,\dots,\phi_n$ as follows. For every $i \in [n]$ we obtain $\phi_i(v_{i-1},v_i)$ from $R_{{i}}(s_i)$ (seen as an atomic formula) by existentially quantifying all variables except for $v_{i-1}$ and $v_i$.   
Furthermore, for every $i\in[0,n]$ and every neighbour $(t,R)$ of $v_i$ in the incidence graph of $\bC$ that is not in $P$, let $\phi(x)$ be the formula
obtained from $R(t)$ by existentially quantifying all variables except for $v_i$.

Let $h$ be a homomorphism from $\bC$ to $\bI$. {Note that the element $h(v_i)$ of the indicator structure $\bI$ is a set of $K$-tuples of elements in $B$. To be able to better work with $h$ we will choose representatives of these sets. We pick them in such a way that they witness that $h$ is a homomorphism. See Figure~\ref{fig:explainTiNotation} for a visualization of the chosen tuples in an example.   
Formally, fix tuples $t_{0}^0,\dots,t_{0}^{m+1},t_{1}^0,\dots,t_{1}^{m+1},\dots,t_{n}^0,\dots,t_{n}^{m+1}$ in $\bB^k$ such that} 
\begin{enumerate}
\item[(P1)] {for every $c \in [0,n]$ we have $t_{c}^0,\dots,t_{c}^{m+1}\in h(v_c)\subseteq B^K$,}
    \item[(P2)] for every $c \in [0,n]$ and every formula $\phi(x)$ corresponding to a neighbour $(t,R)$ of $v_c$ in the incidence graph of $\bC$ that is not in $P$ there is an $i\in[m]$ such that {$\bB^k\models\phi(t_c^i)$} 
    and
    \item[(P3)] {for every $c\in [0,n-1]$ we have  $\bB^k\models\phi_c(t_{c}^{m+1},t_{c+1}^0)$.} 
\end{enumerate}
{
The existence of such tuples follows directly from $h$ being a homomorphism and $\bI$ being an  indicator structure. 
Note that   $\bB^k\models\phi(t_c^i)$ is equivalent to $\bB\models \phi((t_c^i)_k)$ for all $k\in[K]$ and that, similarly,  $\bB^k\models\phi_c(t_{c}^{m+1},t_{c+1}^0)$ is equivalent to $\bB\models \phi_c((t_{c}^{m+1})_k,(t_{c+1}^0)_k)$ for all $k\in[K]$.  
} 

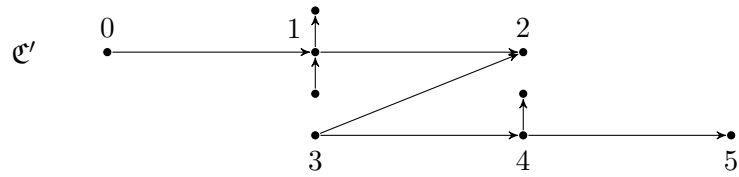
\begin{figure}
    \centering
        \begin{tikzpicture}[scale=0.55]
    
        \node at (0.5,0) {$\bC$};
        \node[var-b, label={below:$v_0$}] (a) at (2.5,0) {};
        \node[var-b, label={below left:$v_1$}] (2) at (7.5,0) {};
        \node[var-b] (1) at (7.5,-1) {};
        \node[var-b] (3) at (7.5,1) {};
        \node[var-b, label={below:$v_2$}] (b) at (12.5,0) {};
        \node[var-b] (b1) at (12.5,1) {};
        \node[var-b, label={below:$v_3$}] (4) at (17.5,0) {};
        
        \path[->,>=stealth']
            (a) edge (2)
            (1) edge (2)
            (2) edge (3)
            (2) edge (b)
            (b) edge (b1)
            (b) edge (4)
            ;
    \end{tikzpicture}\vspace{12mm}
    \begin{tikzpicture}[scale=0.55,rotate=-90]
        \def\colA{black}
        \def\colB{green!85!black}
        \def\colC{orange!90!black}
        \def\colD{blue}
        \def\colE{cyan!80!black}
        \def\colF{red}
        \fill[lightgray] (-.4,-.4) rectangle (2.4,3.4);
        \fill[lightgray] (-.4,-.4+5) rectangle (2.4,3.4+5);
        \fill[lightgray] (-.4,-.4+10) rectangle (2.4,3.4+10);
        \fill[lightgray] (-.4,-.4+15) rectangle (2.4,3.4+15);
        \node[\colA] at (0,0) {0};
        \node[\colB] at (1,0) {3};
        \node[\colB] at (2,0) {3};
        
        \node[\colA] at (0,1) {0};
        \node[\colB] at (1,1) {3};
        \node[\colB] at (2,1) {3};
        
        \node[\colC] at (0,2) {1};
        \node[\colA] at (1,2) {0};
        \node[\colC] at (2,2) {1};
        
        \node[\colC] (11b) at (0,3) {1};
        \node[\colA] (12b) at (1,3) {0};
        \node[\colC] (13b) at (2,3) {1};

        \node[\colC] (21a) at (0,5) {2};
        \node[\colA] (22a) at (1,5) {1};
        \node[\colC] (23a) at (2,5) {2};
        
        \node[\colA] at (0,6) {1};
        \node[\colC] at (1,6) {2};
        \node[\colC] at (2,6) {2};
        \node[rotate=-90] at (4.6,6) {$\in(\exists y E(x,y))^{\bB}$};
        
        \node[\colA] at (0,7) {1};
        \node[\colD] at (1,7) {4};
        \node[\colD] at (2,7) {4};
        \node[rotate=-90] at (4.6,7) {$\in(\exists y E(y,x))^{\bB}$};

        \node[\colA] (21b) at (0,8) {1};
        \node[\colA] (22b) at (1,8) {3};
        \node[\colA] (23b) at (2,8) {3};

        \node[\colA] (31a) at (0,10) {3};
        \node[\colA] (32a) at (1,10) {2};
        \node[\colA] (33a) at (2,10) {4};
        
        \node[\colA] at (0,11) {3};
        \node[\colA] at (1,11) {2};
        \node[\colA] at (2,11) {4};
        \node[rotate=-90] at (4.6,11) {$\in(\exists y E(x,y))^{\bB}$};
        
        \node[\colA] at (0,12) {3};
        \node[\colA] at (1,12) {2};
        \node[\colA] at (2,12) {4};
        
        \node[\colA] (31b) at (0,13) {3};
        \node[\colA] (32b) at (1,13) {2};
        \node[\colA] (33b) at (2,13) {4};

        \node[\colA] (41a) at (0,15) {2};
        \node[\colA] (42a) at (1,15) {2};
        \node[\colA] (43a) at (2,15) {5};
        
        \node[\colE] at (0,16) {1};
        \node[\colA] at (1,16) {5};
        \node[\colE] at (2,16) {1};
        
        \node[\colE] at (0,17) {1};
        \node[\colA] at (1,17) {5};
        \node[\colE] at (2,17) {1};
        
        \node[\colF] (41b) at (0,18) {4};
        \node[\colF] (42b) at (1,18) {4};
        \node[\colA] (43b) at (2,18) {5};

\path[->,>=stealth']
(11b) edge (21a)
(12b) edge (22a)
(13b) edge (23a)
(11b) edge (23a)
(13b) edge (21a)

(21b) edge  (31a)
(22b) edge  (32a)
(23b) edge  (33a)
(21b) edge  (32a)
(22b) edge  (33a)
(23b) edge  (32a)

(31b) edge  (41a)
(32b) edge  (42a)
(33b) edge  (43a)
(31b) edge  (42a)
(32b) edge  (41a)
;
    \end{tikzpicture}\vspace{12mm}
    \begin{tikzpicture}[scale=0.55,rotate=-90]
        
        \fill[lightgray] (-.4,-.4) rectangle (2.4,3.4);
        \fill[lightgray] (-.4,-.4+5) rectangle (2.4,3.4+5);
        \fill[lightgray] (-.4,-.4+10) rectangle (2.4,3.4+10);
        \fill[lightgray] (-.4,-.4+15) rectangle (2.4,3.4+15);

        \def\colA{black}
        \node[\colA] at (0,0) {$a_0$};
        \node[rotate=-90] at (4,0) {$=(0,0,1)$};
        
        \node[\colA] at (0,1) {$a_1$};
        \node[rotate=-90] at (4,1) {$=(0,1,1)$};

        \node[\colA] at (1,2) {$a_2$};
        \node[rotate=-90] at (4,2) {$=(0,2,2)$};

        \node[\colA] (12b) at (1,3) {$a_3$};
        \node[rotate=-90] at (4,3) {$=(0,3,2)$};

        \node[\colA] (22a) at (1,5) {$a_4$};
        \node[rotate=-90] at (4,5) {$=(1,0,2)$};

        \node[\colA] at (0,6) {$a_5$};
        \node[rotate=-90] at (4,6) {$=(1,1,1)$};

        \node[\colA] at (0,7) {$a_6$};
        \node[rotate=-90] at (4,7) {$=(1,2,1)$};

        \node[\colA] (21b) at (0,8) {$a_7$};
        \node[\colA] (22b) at (1,8) {$a_9$};

        \node[\colA] (32a) at (1,10) {$a_8$};
        \node[\colA] (33a) at (2,10) {$a_{10}$};
        
        \node[\colA] at (2,11) {$a_{11}$};
        
        \node[\colA] at (2,12) {$a_{12}$};
        
        \node[\colA] (33b) at (2,13) {$a_{13}$};

        \node[\colA] (43a) at (2,15) {$a_{14}$};
        
        \node[\colA] at (1,16) {$a_{15}$};
        
        \node[\colA] at (1,17) {$a_{16}$};
        
        \node[\colA] (43b) at (2,18) {$a_{17}$};

        \path[->,>=stealth']
(12b) edge (22a)

(22b) edge  (32a)
(21b) edge  (32a)
(22b) edge  (33a)

(33b) edge  (43a)
;
    \end{tikzpicture}\vspace{12mm}
    \begin{tikzpicture}[scale=0.61,label distance=-1.5mm]
        \node[label=right:{$=\{0,1,2,3\}$}] (a) at (2.5,0) {$A_0$};
        \node[label=right:{$=\{4,5,6,7\}$}] (2) at (7.5,0) {$A_1$};
        \node[label=right:{$=\{8\}$}] (b) at (12.5,0) {$A_2$};
        
        \node[label=right:{$=\{9\}$}] (2b) at (7.5,-1) {$A_3$};
        \node[label=right:{$=\{10,\dots,13\}$}] (bb) at (12.5,-1) {$A_4$};
        \node[label=right:{$=\{14,\dots,17\}$}] (4b) at (17.8,-1) {$A_5$};
    \end{tikzpicture}\vspace{12mm}
    \begin{tikzpicture}[scale=0.55]
        \node at (0.5,0) {$\bC'$};
        \node[var-b] (1) at (7.5,-1) {};
        \node[var-b, label={above:$0$}] (a) at (2.5,0) {};
        \node[var-b, label={above left:$1$}] (2) at (7.5,0) {};
        \node[var-b] (3) at (7.5,1) {};
        \node[var-b, label={above:$2$}] (b) at (12.5,0) {};
        
        \node[var-b, label={below:$3$}] (2b) at (7.5,-2) {};
        \node[var-b, label={below:$4$}] (bb) at (12.5,-2) {};
        \node[var-b, label={below:$5$}] (4b) at (17.5,-2) {};
        \node[var-b] (cb) at (12.5,-1) {};

        \path[->,>=stealth']
            (1) edge (2)
            (a) edge (2)
            (2) edge (3)
            (2) edge (b)
            (2b) edge (b)
            (2b) edge (bb)
            (bb) edge (4b)
            (bb) edge (cb)
            ;
    \end{tikzpicture}
    \caption{
    At the top there is a picture of a caterpillar $\bC$.
    Below there is a visualisation of the tuples {$t_0^0,\dots,t_3^3$} introduced in the proof of Lemma~\ref{lem:interesting}. The first column is {$t_0^0$, the second $t_0^1$}, and so on. The condition used to identify tuples is the ternary quasi minority condition. The colors indicate the equivalence classes of the $\connects$ relation.  Below we have the partially labelled sequence $a_0,\dots,a_{17}$ and the resulting sequence $A_0,\dots,A_5$. At the bottom there is the constructed unfolding of $\bC$.}
    \label{fig:explainTiNotation}
\end{figure}

Define binary relations $\sim$ and $\toEdge_a$ on $[0,n]\times [0,m+1]\times [K]$ for $a\in[n]$: 
\begin{itemize}
    \item $(c,i,k)\sim (c,i+1,\ell)$ if and only if {$(t_c^i)_k=(t_c^{i+1})_\ell$, note that the tuples $t_c^i$ and $t_c^{i+1}$ represent the same element in $\bI$}
    \item $(c,m+1,k) \toEdge_a (c+1,0,\ell)$  if and only if  {$\bB\models\phi_a((t_c^{m+1})_k,(t_{c+1}^0)_\ell)$}, {note that because of (P3) this conditions holds whenever $k=\ell$, though it may also hold if $k\neq\ell$.}
\end{itemize}
We write $\fromEdge_a$ for the converse of $\toEdge_a$. Define \connects\ as the {reflexive} transitive symmetric closure of $\sim\cup\toEdge_1\cup\dots\cup\toEdge_n$. Note that \connects\ is an equivalence relation. 
Now we consider two cases. The first case
is that there are no $k,\ell\in[K]$ such that $(0,0,k) \connects (n,m+1,\ell)$. Define the map 
\begin{align*}
    \pi  \colon [0,n]\times [0,m+1]\times [K] &\to\{x,y\}\\
    (c,i,k)&\mapsto 
    \begin{cases}
        x&\text{if $(c,i,k) \connects (0,0,\ell)$ for some $\ell$}\\
        y&\text{otherwise.}
    \end{cases}
\end{align*}   
Observe that 
\begin{enumerate}
    \item {By assumption} $\pi((0,0,k))=x$ for all $k\in[K]$ and that $\pi((n,m+1,\ell))=y$ for all $\ell\in[K]$.
    \item Let $c\in[0,n]$ and let $i\in[0,m]$. By (P1) $t_c^i$ and $t_c^{i+1}$ represent the same element in the indicator structure $\bI$. Hence, when viewing the entries of $t_c^i$ and $t_c^{i+1}$ as variables, $\Sigma$ implies $f(t_c^i)\approx f(t_c^{i+1})$. By definition of $\pi$ and $\sim$ we have that $(t_c^{i})_{k_1}=(t_c^{i+1})_{k_2}$ implies $\pi(c,i,{k_1})=\pi(c,i+1,{k_2})$ for all $k_1,k_2$. Therefore the minor condition $f(t_{c}^i)\approx f(t_{c}^{i+1})$ implies the minor condition 
    \[f(\pi(c,i,1),\dots,\pi(c,i,K))\approx f(\pi(c,i+1,1),\dots,\pi(c,i+1,K)).\]
    \item Let $c\in[0,n-1]$, then, {by definition of $\toEdge_a$ and (P3), we have that}  $\pi(c,m+1,k)=\pi(c+1,0,k)$ for all $k\in[K]$. Hence, 
    \[f(\pi(c,m+1,1),\dots,\pi(c,m+1,K))= f(\pi(c+1,0,1),\dots,\pi(c+1,0,K)).\]
\end{enumerate}
Therefore, $\Sigma$ implies 
\begin{align*}
f(x,\dots,x) & =f(\pi(0,0,1),\dots,\pi(0,0,K))\\
& \approx f(\pi(0,1,1),\dots,\pi(0,1,K))\\
& \quad \vdots\\
& \approx f(\pi(n,m+1,1),\dots,\pi(n,m+1,K))=f(y,\dots,y).
\end{align*}
Hence, $\Sigma$ implies $g(x)\approx g(y)$, as desired.

Now we consider the case that there are $k,\ell\in[K]$ such that $(0,0,k) \connects (n,m+1,\ell)$. 
Then there exists a sequence  $a_0,\dots,a_N$ of distinct elements in $[0,n]\times[0,m+1]\times[K]$ such that $a_0=(0,0,k)$, $a_N=(n,m+1,\ell)$, and for every $i\in[N]$ we have $a_{i-1}\sim a_i$, $a_{i-1}\toEdge_j a_i$, or $a_{i-1}\fromEdge_j a_i$ for some $j\in[n]$.  
We will use this sequence to construct an unfolding of $\bC$ that has a homomorphism to $\bB${, contradicting that $\bB$ has unfolded caterpillar duality while $\bC\not\to\bB$.} See Figure~\ref{fig:explainTiNotation} 
for an example of such a sequence and the construction that follows.  
First group the sequence $0,\dots,N$ into the sequence $A_0,\dots,A_M$ such that 
\begin{itemize}
    \item $a_0\in A_0$ and $a_N\in A_M$,
    \item {each $A_i$ consists of an interval, i.e., for all $i\in[0,M]$ there are $0\leq j_1\leq j_2\leq N$ such that $A_i=\{{j_1},{j_1+1},\dots,{j_2}\}$},
    \item for all $i\in[0,M]$ and $j\in A_i$  with $j+1\in A_i$ we have $a_{j}\sim a_{j+1}$,
    \item for all $0\leq i_1<i_2\leq M$ and all ${j_1}\in A_{i_1}$, ${j_2}\in A_{i_2}$ we have $j_1<j_2$, and
    \item for all $i\in[M]$ we have $\max(A_{i-1})+1=\min(A_i)$ and there is a $j\in[n]$ with $a_{\max(A_{i-1})}\toEdge_j a_{\min(A_i)}$ or $a_{\max(A_{i-1})}\fromEdge_j a_{\min(A_i)}$.
\end{itemize}
We write $A_{i}\toEdge_j A_{i+1}$ if $a_{\max(A_{i-1})}\toEdge_j a_{\min(A_i)}$ and  $A_{i}\fromEdge_j A_{i+1}$ if $a_{\max(A_{i-1})}\fromEdge_j a_{\min(A_i)}$. 
Define the primitive positive formula $\Phi(0,\dots,M)$ by adding conjuncts in the following way: 
For all $i$ and all $j$
\begin{itemize}
    \item if $A_{i}\toEdge_j A_{i+1}$, then add the conjunct $\phi_j(i,{i+1})$,
    \item if $A_{i}\fromEdge_j A_{i+1}$, then add the conjunct $\phi_j({i+1},{i})$,
    \item if $A_{i-1}\toEdge_j A_{i}\toEdge_{j+1} A_{i+1}$, then for all $\phi(x)$ corresponding to a neighbour $(t,R)$ of $v_j$ in the incidence graph of $\bC$ that is not in $P$ add the conjunct $\phi({i})$, and
    \item if $A_{i-1}\fromEdge_{j+1} A_{i}\fromEdge_{j} A_{i+1}$, then for all $\phi(x)$ corresponding to a neighbour $(t,R)$ of $v_j$ in the incidence graph of $\bC$ that is not in $P$ add the conjunct $\phi({i})$.
\end{itemize}
{Note that in the last two cases, by definition of $\sim$,  \[[0,m+1]=\{i'\mid (c,i',k) \in \{a_j\mid j \in A_i\}\}.\] }
Denote the canonical database of $\Phi$ by $\bC'$. Note that $\bC'$ is a caterpillar.
Observe that, by definition, {for all $i\in[0,M]$ and $j_1,j_2\in A_{i}$ with $a_{j_1}=(c,i',k)$ and $a_{j_2}=(d,j',\ell)$ we have that $t_c^{i'}$ and $t_d^{j'}$ represent the same element in $\bC$ and that $(t_c^{i'})_k=(t_d^{j'})_{\ell}$. Hence, the map $\iota\colon [0,M]\to B$ that maps $i$ to $(t_c^{i'})_k$ for some (any) $(c,i',k)$ for which there is a $j\in A_{i}$ with $(c,i',k)=a_{j}$ is well defined.}   
The map $\iota$ is a satisfying assignment of $\Phi$ in $\bB$:
\begin{itemize}
    \item {If $\Phi$ contains the conjunct $\phi_j(i,i+1)$, then $A_i\toEdge_j A_{i+1}$. Hence $(c,i',k)=a_{\max(A_i)}\toEdge_j a_{\min(A_{i+1})}=(d,j',\ell)$ and $\iota(i)=(t_c^{i'})_k$ and $\iota(i+1)=(t_d^{j'})_{\ell})$. By definition of $\toEdge_j$ we have $\bB\models\phi_j((t_c^{i'})_k,(t_d^{j'})_{\ell})$. }
    \item The argument for conjuncts of the form $\phi_j(i+1,i)$ is analogous.
    \item {The conjuncts of the form $\phi(i)$ are satisfied: Let $(d,j,\ell)=a_{\min(A_i)}$. Then, by (P2), there is a $j'\in[m]$ such that $\bB\models \phi((t_d^{j'})_k)$ for all $k\in[K]$. By definition of $\Phi$ we added $\phi(i)$ only if there is some $j\in A_i$ with $a_{j}=(d,j',k)$ for some $k\in[K]$.
    }
    
\end{itemize}
Therefore, $\iota$ can be extended to a homomorphism from $\bC'$ to $\bB$.
The proof that $\bC'$ is an unfolding of $\bC$ is technical and can for example be done by induction on the number of \emph{orientation changes} in the sequence $A_0,\dots, A_M$. However, it is easy to see that $\bC'$ must indeed be an unfolding of $\bC$. 
\end{proof}

\begin{remark}
    The conclusion of Lemma~\ref{lem:interesting} can be strengthened as follows: if $\bB$ does not satisfy a minor condition $\Sigma$, then $\Sigma$ implies $f(x) \approx f(y)$ 
    even with respect to the class of all minions, i.e., it then holds that every \emph{minion} that satisfies $\Sigma$ also satisfies $f(x) \approx f(y)$. However, in this case we cannot use Lemma~\ref{lem:one-f} to assume that $\Sigma$ uses only a single function symbol, since it is a statement for clones rather than minions. 
    Instead, one can then use a more general notion of indicator structure and adapt the entire proof to this more general setting. Since we do not need this for proving our main result,  Theorem~\ref{thm:main}, we have decided for the weaker result which allows for a less technical proof. 
\end{remark}

\subsection{Proof of the Main Result}
We finally prove Theorem~\ref{thm:main}. 

\begin{proof}[Proof of Theorem~\ref{thm:main}]


For the implication $(\ref{maltsev})\Rightarrow(\ref{can-sym-lin-arc})$ let $\Pi_S$ be the canonical slam Datalog program for $\bB$ and let $\Pi$ be the canonical linear arc monadic Datalog program for $\bB$. Since $\bB$ has $k$-absorptive operations of arity $nk$ for all $n,k\geq 1$ we can apply  Theorem~\ref{thm:caterpillar} to conclude that $\Pi$ solves $\Csp(\bB)$. Furthermore, $\bB$ has a quasi Maltsev polymorphism. Hence, Lemma~\ref{lem:MaltimpliesLAMDatalogIsSLAMDatalog} implies that $\Pi$ and $\Pi_S$ can derive the goal predicate on the same instances of $\Csp(\bB)$. Therefore, $\Pi_S$ solves $\Csp(\bB)$.


The implication $(\ref{can-sym-lin-arc})\Rightarrow(\ref{sym-lin-arc})$ is trivial, the implication 
$(\ref{sym-lin-arc})\Rightarrow(\ref{caterpillar})$
by Lemma~\ref{lem:SLAMimpliesUnfoldedCaterpillar}, and the implication  $(\ref{caterpillar})\Rightarrow(\ref{non-degenerate-minor})$ 
 by Lemma~\ref{lem:interesting}.


For the implication from~\eqref{non-degenerate-minor} to~\eqref{minor-p2}, suppose that $\Sigma$ is a minor condition that holds in $\Pol(\bP_2)$. 
Since all polymorphisms of $\bP_2$ are idempotent, 
$\Sigma$ does not imply $f(x) \approx f(y)$. 
Hence, the contraposition of~\eqref{non-degenerate-minor} implies that $\Pol(\bB)$ does not satisfy $\Sigma$. 


The implication from~\eqref{minor-p2} to~\eqref{maltsev} is clear since $\Pol(\bP_2)$ is preserved by the Boolean minority operation and by the $nk$-ary Boolean operation
$$(x_{11},\dots,x_{1k},\dots,x_{n1},\dots,x_{nk}) \mapsto \bigvee_{i \in [n]} \bigwedge_{j \in [k]} x_{ij}$$
which is $k$-absorptive. 


The equivalence between~\eqref{minor-p2}, \eqref{minion-hom}, and~\eqref{pp-p2} follows from Remark~\ref{rem:p2} and  well-known general results~\cite{wonderland}.


The equivalence of~\eqref{maltsev} and~\eqref{lattice} follows the equivalence of items (3) and (4) in Theorem~\ref{thm:caterpillar} and from the fact that the existence of a quasi Maltsev polymorphism is preserved by homomorphic equivalence.
\medskip 

For the final statement of the theorem, {assume that $\bB$ is a core.} 
Let $\bB'$ be the structure with the same domain as $\bB$ which contains all binary relations that are primitively positively definable in $\bB$. 
First recall from Remark~\ref{rem:majo} that \eqref{maltsev} implies that $\bB$ has a quasi majority polymorphism, and hence every relation of $\bB$ is equivalent to a conjunction of binary relations of 
$\bB'$ (Proposition~\ref{prop:decomp}), which shows that $\Pol(\bB) = \Pol(\bB')$. 
In particular, $\bB'$ has $k$-absorptive polymorphisms 
        of arity $nk$, for all $n,k \geq 1$, and hence the theorem applies to $\bB'$ in place of $\bB$ as well. 
\end{proof}

\begin{remark}\label{rem:T3}
    Consider the poset of all finite structures ordered by primitive positive constructability. It is well known that the structure $\bC_1\coloneqq (\{0\},\{(0,0)\})$ is a representative of the top element of this poset and that it has exactly one lower cover with representative {$\bP_2$}. 
    We claim that $\bT_3$ is a representative of a lower cover of $\bP_2$ in the poset of all finite structures ordered by primitive positive constructability. 
    The structure $\bT_3\coloneqq(\{0,1,2\},\{(0,1),(0,2),(1,2)\})$ satisfies all conditions $\Sigma$ that do not imply the quasi Maltsev condition (see, e.g.,~\cite{maximal-digraphs}). 
Clearly, $\bT_3$ does not have a primitive positive construction in $\bP_2$, because $\bT_3$ does not have a quasi Maltsev polymorphism. 
    Since $\min$ and $\max$ are polymorphisms of $\bT_3$, Theorem~\ref{thm:caterpillar} implies that $\bT_3$ has $kn$-ary $k$-absorbing polymorphisms for all $n,k\geq 1$. 
    Let $\bB$ be a structure with a primitive positive construction in $\bT_3$ which does not admit a primitive positive construction of $\bT_3$. 
    Then $\bB$ must have a quasi Maltsev polymorphism and $kn$-ary $k$-absorbing polymorphisms for all $n,k\geq 1$. By Theorem~\ref{thm:main}, $\bB$ has a primitive positive construction in $\bP_2$, which proves the claim.
    It is still open what other lower covers $\bP_2$ has.
\end{remark}

\section{Decidability of Meta-Problem}
There are many interesting results and open problems about \emph{algorithmic meta-problems} in constraint satisfaction; we refer to~\cite{MetaChenLarose}. The natural algorithmic meta-problem in the context of our work is the one addressed in the following proposition. 


\begin{proposition}\label{prop:meta}
    There is an algorithm which decides in deterministic doubly-exponential time whether the CSP of a given finite structure $\bB$ can be solved by a slam Datalog program, and if so, computes such a program. 
\end{proposition}
\begin{proof}
The following algorithm can be used to test whether $\bB$ has $k$-absorptive operations of arity $nk$, for all $n,k \geq 1$. It is well-known that the existence of a quasi Maltsev polymorphism can be decided in non-deterministic polynomial time (see, e.g.,~\cite{MetaChenLarose}). 
Hence, the statement then follows from 
Theorem~\ref{thm:main}.

Let $m$ be the maximal arity of the relations of $\bB$.
Let $n_0 := m\binom{|B|}{|B|/2}$ and $k_0 := m|B|$.
Note that $\bB$ has $k$-absorptive polymorphisms of arity $nk$, for all $k,n$, if and only if it has 
$k_0$-absorptive polymorphisms of arity $n_0 k_0$ (similarly as the well-known fact that $\bB$ has totally symmetric polymorphisms of all arities if and only if it has totally symmetric polymorphisms of arity $m|B|$; the term $\binom{|B|}{|B|/2}$ bounds the size of antichains in the set of all subsets of $B$. Also see~\cite{CarvalhoDalmauKrokhin} for the case of $k$-absorptive polymorphisms). Let $\Sigma$ be the minor condition for the existence of $k_0$-absorptive operations of arity $n_0k_0$. 
Let $\bC$ be the indicator structure 
of $\Sigma$ with respect to $\bB$ as defined in Section~\ref{sect:indicator};
clearly, this structure can be computed in doubly exponential time. 
We may then find a non-deterministic algorithm with the same time bound that tests whether there exists a homomorphism from $\bC$ to $\bB$. 
 The non-determinism for checking whether $\bC \to \bB$ 
can be eliminated by standard self-reduction techniques (again see, e.g.,~\cite{MetaChenLarose}). 
For the second part of the statement, note 
that there are for a given $\bB$ only finitely many potential rules of a slam Datalog program, and one can compute for a given rule whether it is part of the canonical slam Datalog program of $\Csp(\bB)$.
\end{proof}

\section{Remarks on Related Results} 
The following remarks show that the results  
of Carvalho, Dalmau and Krokhin~\cite{CarvalhoDalmauKrokhin} can be extended in the same spirit as our Theorem~\ref{thm:main}. 
    
\begin{remark}\label{rem:caterpillar}
Let $\bD_2$ be the structure $(\{0,1\};\{0\},\{1\},\leq)$, {where $\leq$ denotes the binary relation $\{(0,0),(0,1),(1,1)\}$}, also known as st-Con. 
Theorem~\ref{thm:caterpillar} 
    of Carvalho, Dalmau and Krokhin can be extended in the same spirit as our Theorem~\ref{thm:main}, by adding the following equivalent items:
    \begin{enumerate}
        \setcounter{enumi}{4}
        \item \label{item:minor-cond-d2} Every minor condition that holds in 
        $\Pol(\bD_2)$ also holds in $\Pol(\bB)$. 
        \item \label{item:minion-homo-d2} There is a minion homomorphism from $\Pol(\bD_2)$ to $\Pol(\bB)$.
            \item \label{item:pp-constr-d2} $\bB$ has a primitive positive construction in $\bD_2$. 
    \end{enumerate}
    The equivalence of \ref{item:minor-cond-d2}., \ref{item:minion-homo-d2}., and \ref{item:pp-constr-d2}.\  follows immediately from the general results in~\cite{wonderland}. 
    
    $\ref{item:cp-4}. \Rightarrow \ref{item:pp-constr-d2}.$ 
    It is well known that 
    $\Pol(\bD_2)$ is generated by the two binary operations $\vee$ and $\wedge$.\footnote{Proof sketch: clearly, $\vee$ and $\wedge$ preserve the relations of $\bD_2$. For the converse inclusion, 
    {note that every Boolean relation preserved by 
    $\vee$ and $\wedge$ has a definition in CNF which is both Horn and dual Horn, so consists of clauses that can be defined using the relations in $\bD_2$ (see, e.g.~\cite{KaluzninPoeschel})}.}
    Let $\bB$ be a structure that is homomorphically equivalent to a structure $\bB'$ with binary polymorphisms $\sqcup$ and $\sqcap$ such that $(B',\sqcup,\sqcap)$ is a distributive lattice. 
    Note that $(\{0,1\},\vee,\wedge)$ is a distributive lattice as well. Let $\iota$ be the map that maps terms over $\vee,\wedge$ to terms over $\sqcup,\sqcap$ by replacing $\vee$ and $\wedge$ by $\sqcup$ and $\sqcap$, respectively. 
    Define the map $\xi\colon \Pol(\bD_2)\to\Pol(\bB')$ as follows. 
    Since $\Pol(\bD_2)$ is generated by $\vee$ and $\wedge$, for every 
    $f\in \Pol(\bD_2)$
    there is a $\{\wedge,\vee\}$-term $t$ whose term operation is $f$. Define $\xi(f)$ as the term operation of $\iota(t)$. Note that this term operation is a polymorphism of $\bB'$. 
    It is clear that $\xi$ is a minion homomorphism (even a clone homomorphism). We still need to show that $\xi$ is well defined. Let $t$ and $t'$ be two $\{\wedge,\vee\}$-terms that both have the term operation $f\in\Pol(\bD_2)$. Since $(\{0,1\},\vee,\wedge)$ is a distributive lattice, there is a set $\mathcal I$ of subsets of $[n]$ such that $f$ is the term operation of 
    \begin{align*}
    s\coloneqq \bigwedge_{I\in\mathcal I}\bigvee_{i\in I} x_i.
    \end{align*}
    Furthermore, $t$ and $t'$ can both be rewritten (using associativity, commutativity, distributivity, and idempotence) into the term $s$. Therefore, $\iota(t)$ and $\iota(t')$ can also both be rewritten into the term $\iota(s)$. Since $(B',\sqcup,\sqcap)$ is a distribute lattice, the 
    term operations of $\iota(t)$, $\iota(t')$, and $\iota(s)$ are the same. Hence, $\xi$ is well defined.


    $\ref{item:minor-cond-d2}. \Rightarrow \ref{item:cp-4}.$ holds since $\bD_2$ has for every $n,k\geq1$ a $k$-absorbing polymorphism of arity $kn$.



\end{remark}

\begin{remark}
Let $\bB_{\infty}^\leq$ be the relational structure with  the domain $\{0,1\}$ and the signature $\{{\boldsymbol 0},\leq, R_1,R_2,\dots\}$ where
 ${\boldsymbol 0}\coloneqq\{0\}$,
 $\leq$ is as in Remark~\ref{rem:caterpillar}, and
    $R_n\coloneqq \{0,1\}^n\setminus \{(0,\dots,0)\}$ for every $n\geq1$. 
It is well known that  $\Pol(\bB_{\infty}^\leq)$ is generated by the operation $m$ given by $(x,y,z)\mapsto x\wedge(y\vee z)$.\footnote{Proof sketch: clearly, every relation of $\bB_{\infty}^\leq$ is preserved by $m$. For the converse inclusion, it suffices to verify that every relation that is preserved by $m$ has a primitive positive definition in $\bB^{\leq}_{\infty}$ (see, e.g.,~\cite{KaluzninPoeschel}). 
First note that $m(x,y,y) = x \wedge y$, and hence every Boolean relation $R$ preserved by $m$ has a Horn definition; pick such a definition $\phi$ which is shortest possible. Suppose for contradiction that a Horn clause in $\phi$ contains a positive literal $\psi_1$ and two negative literals $\psi_2$ and $\psi_3$. By the minimality assumption there are tuples $t_1,t_2,t_3 \in R$ such that $t_i$ satisfies $\phi_i$ and no other literal in that clause. Then $m(t_1,t_2,t_3)$ satisfies none of  $\psi_1,\psi_2,\psi_3$, a contradiction. It follows that each clause can be defined using the relations in $\bB_{\infty}^\leq$ and the statement follows.}
Carvalho, Dalmau and Krokhin also introduce another type of duality in their paper: jellyfish duality. Their characterization in Theorem~18 in~\cite{CarvalhoDalmauKrokhin} can be extended by the following items:
    \begin{enumerate}
        \setcounter{enumi}{5}
        \item Every minor condition that holds in 
        $\Pol(\bB_{\infty}^\leq)$ also holds in $\Pol(\bB)$. 
        \item There is a minion homomorphism from 
        $\Pol(\bB_{\infty}^\leq)$ to $\Pol(\bB)$.
            \item $\bB$ has a primitive positive construction in $\bB_{\infty}^\leq$. 
    \end{enumerate}
    The proof is analogous to the proof in Remark~\ref{rem:caterpillar}.
\end{remark}

\begin{figure}
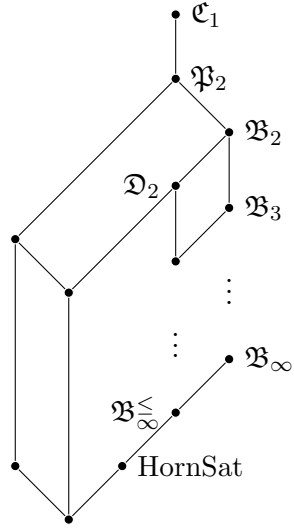

    \centering
    \albertTwoElementPoset
    \caption{The lattice of 2-element structures with respect to pp-constructability~\cite{PPPoset}; the structures HornSat, $\bD_2$, and $\bP_2$ were relevant in this text for characterisations of Datalog fragments (for arc monadic, linear arc monadic, and slam Datalog, respectively).}
    \label{fig:twoElementPoset}
\end{figure}

\begin{remark}\label{rem:gm}
    Another condition on finite structures $\bB$ that is equivalent to the conditions in Theorem~\ref{thm:main} has been found by Vucaj and Zhuk~\cite{vucaj2024submaximal}: they prove that there is a minion homomorphism from $\Pol(\bP_2)$ to $\Pol(\bB)$ (item~\ref{minion-hom}) if and only if $\Pol(\bB)$ contains totally symmetric polymorphisms of all arities and \emph{generalised quasi minority polymorphisms} of all odd arities $n \geq 3$. 
    A operation $f \colon B^n \to B$, for odd $n \geq 3$, is called a \emph{generalised quasi minority} if
    it satisfies 
    \begin{align*}
        f(x_1,\dots,x_n) & \approx f(x_{\pi(1)},\dots,x_{\pi(n)}) && \text{ for every } \pi \in S_n \\
        \text{ and } f(x,x,x_3,\dots,x_n) & \approx f(y,y,x_3,x_4,\dots,x_n) . 
    \end{align*}
    However, it is not clear to us whether this characterisation can be used to prove the consequence of our main result from Remark~\ref{rem:T3}.
\end{remark}

\begin{remark}
Yet another remarkable equivalent condition 
for solvability by slam Datalog was very recently discovered by Meyer and Starke~\cite{meyerStarke2024finitesimplegroupsprimitive}: the conditions from Theorem~\ref{thm:main} hold if and only if 
    $\bB$ can neither pp-construct $\bT_3$ nor any structure from a list of structures that  correspond to the finite simple groups with a particular action; for details, we refer to~\cite{meyerStarke2024finitesimplegroupsprimitive}.  
\end{remark}

\begin{remark}
A minion $\mathscr M$ is called a \emph{core} if every minion homomorphism from $\mathscr M$ to $\mathscr M$ 
is injective. 
We say that $\mathscr N$ is a \emph{minion core of $\mathscr M$} if 
$\mathscr M$ and $\mathscr N$ are homomorphically equivalent (i.e., there is a homomorphism from $\mathscr M$ to $\mathscr N$ and vice versa)
and $\mathscr N$ is a minion core. 
Note that if $\mathscr M$ is \emph{locally finite}, i.e., 
if $\mathscr M^{(n)}$ is finite for every $n \in {\mathbb N}$, then
there exists a minion $\mathscr N$ 
which is a minion core of $\mathscr M$, and $\mathscr N$ is unique up to isomorphism. We therefore call it \emph{the} minion core of $\mathscr M$. 
From personal communication with Libor Barto we learned that the equivalent items of Theorem~\ref{thm:main} apply if and only if the minion core of $\Pol(\bB)$ equals $\Pol(\bC_2,\bB_2)${, where 
\begin{itemize} 
\item $\bC_2$ is the structure $(\{0,1\},\{1\},\{(0,1),(1,0)\})$,
\item $\bB_2$ is the structure $(\{0,1\},\{1\},\{(0,0),(0,1),(1,0)\})$, and 
\item $\Pol(\bA,\bB)$, for relational structures $\bA$ and $\bB$ with the same signature, is the set of all homomorphisms from $\bA^k$ to $\bB$, for all $k \in {\mathbb N}$.
\end{itemize}}
\end{remark}

\section{Conclusion and Open Problems}We characterised the unique submaximal element 
in the primitive positive constructability poset on finite structures {(represented by $\bP_2$)}, linking concepts from homomorphism dualities, Datalog fragments, minor conditions, and minion homomorphisms.
It is now tempting to further descend
in the pp-constructability poset of finite structures
in order to obtain a more systematic understanding. Particularly attractive are other dividing lines in the poset that are relevant for the complexity of the constraint satisfaction problem. We propose the following problems for future research.
\begin{itemize}
    \item Is there a countable set of structures $\bC_1,\bC_2,\dots$ such that $\bB$ does not have a pp-construction from
    $\bP_2$ if and only if one of the structures $\bC_1,\bC_2,\dots$ has a pp-construction from $\bB$? This is true if we restrict to digraphs~\cite{maximal-digraphs} and if we restrict to 3-element structures~\cite{vucaj2024submaximal}. 
    Our result shows that $\bT_3$ must belong to this set (Remark~\ref{rem:T3}). 
    \item Characterise all finite structures that are primitively positively constructible in a finite structure that has finite duality. Are  
    these exactly the finite structures whose polymorphism clones have Hagemann-Mitschke chains of some length 
    and \emph{extended $k$-absorptive polymorphisms of arity $kn+1$}, for all $n,k \geq 1$, as defined in~\cite{Lattice-Ops}? 
   Is there a Datalog fragment that corresponds to this class? 
\item What is the precise computational complexity the Meta-Problem of deciding whether the CSP of a given finite structure $\bB$ can be solved by a slam Datalog program? The algorithm from Proposition~\ref{prop:meta} only provides a deterministic doubly exponential time algorithm.
\end{itemize}

\bibliographystyle{abbrv}
\bibliography{global}

\end{document}